\documentclass{amsart}
\usepackage{amsmath,amssymb,amsfonts,caption}
\usepackage{graphicx}
\usepackage{pstricks}
\usepackage{euscript}
\usepackage{enumerate}
\usepackage{verbatim}
\usepackage{epsfig}
\usepackage{booktabs}
\usepackage{subcaption}

\newtheorem{theorem}{Theorem}
\newtheorem{lemma}{Lemma}

\newtheorem{remark}{Remark}

\newcommand{\eps}{\varepsilon}
\newcommand{\lbd}{\lambda}
\renewcommand{\phi}{\varphi}

\DeclareMathOperator{\e}{e}

\renewcommand{\le}{\leqslant}
\renewcommand{\ge}{\geqslant}
\usepackage{dsfont}
   
   \newcommand{\R}{\ensuremath{\mathds R}}
\begin{document}
\title[A mathematical model for the customer dynamics...]
   {A mathematical model for the customer dynamics based on marketing policy}
\author{C\'esar M. Silva}
\address{
   C\'esar M. Silva\\
   Departamento de Matem\'atica and CMA-UBI\\
   Universidade da Beira Interior\\
   6201-001 Covilh\~a\\
   Portugal}
\email{csilva@ubi.pt}
\author{Silv\'erio Rosa}
\address{
   Silv\'erio Rosa\\
   Departamento de Matem\'atica and IT\\
   Universidade da Beira Interior\\
   6201-001 Covilh\~a\\
   Portugal}
\email{rosa@ubi.pt}
\author{Helena Alves}
\address{
   Helena Alves\\
   Departamento de Gest\~ao e Economia and NECE-UBI\\
   Universidade da Beira Interior\\
   6201-001 Covilh\~a\\
   Portugal}
\email{halves@ubi.pt}
\author{Pedro G. Carvalho}
\address{
   Pedro G. Carvalho\\
   CIDESD-UBI\\
   Universidade da Beira Interior\\
   6201-001 Covilh\~a\\
   Portugal}
\email{pguedes@ubi.pt}
\thanks{C\'esar M. Silva was partially supported by FCT through CMA-UBI (project PEst-OE/MAT/UI0212/2014)}
\begin{abstract}
We consider a compartmental model to study the evolution of the number of regular customers and referral customers in some corporation. Transitions between compartments are modeled by parameters depending on the social network and the marketing policy of the corporation. We obtain some results on the asymptotic number of regular customers and referral customers in several particular scenarios. Additionally we present some simulation that illustrates the behavior of the model and discuss its applicability.
\end{abstract}
\maketitle
\section{Introduction}

In marketing literature it has been successively referred the importance of calculating the value of a customer. In fact, such indicative value enables firms to select those customers that can add profit and consequently constitutes an important information to segment the market and efficiently allocate marketing policy resources~\cite{Kumar-2006a,Kumar-2010,Kumar-Petersen-Leone-2007}.

The objective of this work is to establish and study a compartmental model, mathematically translated into a system of ordinary differential equations,  for the evolution of the number of customers of some firm, assuming that the customers are divided in two subgroups corresponding to different profitabilities.

Until recently, the value of a customer for a company was based on the present value of future profits generated by a customer over the full course of their dealings with a particular company, this is the customer life-time-value (CLV)~\cite{Kumar-al-2010}. However, other authors refer to the importance of including not only the present and future revenue from the customer purchases, but also the value of the potential to influence other customers under incentives on behalf of the company (customer referral value) or by own initiative (customer influencer value)~\cite{Kumar-2010}. Customer influencing behaviors consists of the intrinsic behaviors motivating the customer to persuade and influence other customers without there being any type of reward on behalf of the company and thus designated the customer influencer value (CIV). In turn, the patterns of customer recommendation are related to the acquisition of new customers due to company initiatives that reward recommendations made to other customers, and thereby establishing the customer referral value (CRV).
According to Kumar et al.~\cite{Kumar-al-2010}, these components are mutually interwoven. Thus, CLV positively correlates with CRV (although only up to a certain point and in an inverted U-shaped relational curve, which means customers reporting average CLV are those most interested in company referral programs) and CLV is positively related with CIV (with an inverted U-shaped relationship in effect between these two concepts).
Much of the literature has focused on the customer referral value through the influence customers might have on the formation of other customers’ attitudes (\cite{Bone-1995}) in the purchasing making decision~\cite{Bansal-Voyer-JSR-2000} and in the reduction of other customers perceived risks~\cite{Godes-Mayzlin-MS-2004}, but little is known about how this processes occur. Since customer referral value and influencer value might have a great impact for companies, these latter try to identify the most influential customers~\cite{Kiss-Bichler-DSS-2008}.

A number of studies allow us to think that the customers of a firm can be classified into several groups according to their influential role over other potential buyers. In imperfect competitive markets information is not purely transparent; some persons are more able than others of influencing people to become a customer of that firm. It is also acceptable to assume that knowing the referrals among each firm’s customers and quantifying their influence constitutes an important asset for the firm competitive advantage, although all customers are important, referrals would be more valuable.

We mostly agree with Marti and Zenou~\cite{Marti-Zenou-2009} when they state that physics/applied mathematics are capable of reproducing many real networks but never reach to explain why they emerge; the economists are very precise to explain why they emerge but their approach does a poor job in matching real world networks. That is why some game theorists are now improving models which take networks as given entities and study the impact of their structure on individuals’ outcomes.

Based on the network theory some models have been tested to study the way influential customers can influence other consumers. For instance Kiss and Bichler~\cite{Kiss-Bichler-DSS-2008} tested real network models, simulated networks and diffusion models to predict influence between customers based on their position within the network. However, as the authors mention this analysis not always is possible if we do not know or do not have information regarding the customer social network. Therefore, other models are needed to try to explain these processes.

In this work we propose a model suitable to describe the dynamics of the number of customers of a given firm. This model is given by a system of ordinary differential equations whose variables correspond to groups of customers and potential customers divided according to their profile and whose parameters reflect the structure of the underlying social network and the marketing policy of the firm. We intend to understand the flows between these groups and its consequences on the raise of customers of the firm. We also want to highlight the usefulness of these models in helping firms deciding their marketing policy.

Specifically, the main objectives of our study is threefold: we intend to obtain theoretical results concerning the long term behavior of the number of customers in various scenarios, we want to present some simulation aimed at illustrating the possibilities of application of our model and, finally, we want to discuss the benefits and limitations of this type of analysis.

As referred, we will consider a compartmental model. As far as we are aware, this is the first time such type of mathematical model is considered in the context of marketing research. We believe that this type of model can be fruitfully explored in this context. This believe is based on the fact that compartmental models have proved to be an important tool not only in the natural sciences, particularly in mathematical epidemiology~\cite{Brauer-Driessche-Jianhong-LNM-2008} and in population biology~\cite{Thieme-PUP-2003, Zhao-SV-2003}, but also, with increasing notoriety in recent years, in the context of economy and other social sciences~\cite{Lin-NAHS-2008, Stiglitz-1993, Tramontana-IRE-2010, Artzrouni-Tramontana-JPM-2014}.

\begin{figure}[h]
\begin{picture}(360,190)(-30,-20)
\setlength{\unitlength}{.3mm}
\put(5,0){\framebox(60,45){$C$}}
\put(5,100){\framebox(60,45){$R$}}
\put(270,0){\framebox(60,45){$P_C$}}
\put(270,100){\framebox(60,45){$P_R$}}
\put(80,140){\vector(1,0){180}}
\put(160,145){$\beta_2 R$}
\put(260,120){\vector(-1,0){180}}
\put(160,127){$\lambda_3 P_R$}
\put(260,110){\vector(-1,0){180}}
\put(110,92){$\left( \lambda_2  + m_R \lambda_6 \right) R P_R + m \lambda_4 P_R $}
\put(80,40){\vector(1,0){180}}
\put(160,45){$\beta_1 C$}
\put(260,20){\vector(-1,0){180}}
\put(160,27){$\lambda_1 P_C$}
\put(260,10){\vector(-1,0){180}}
\put(130,-7){$\lambda_2 R P_C + m \lambda_4 P_C$}
\put(25,50){\vector(0,1){45}}
\put(0,70){$\lambda_5 C$}
\put(35,95){\vector(0,-1){45}}
\put(40,70){$\lambda_7 R$}
\put(295,50){\vector(0,1){45}}
\put(265,70){$\lambda_5 P_C$}
\put(305,95){\vector(0,-1){45}}
\put(310,70){$\lambda_7 P_R$}
\put(0,25){\vector(-1,0){30}}
\put(-20,10){$\epsilon C$}
\put(0,115){\vector(-1,0){30}}
\put(-20,100){$\epsilon R$}
\put(375,30){\vector(-1,0){40}}
\put(335,37){$(1-\alpha) \gamma$}
\put(375,130){\vector(-1,0){40}}
\put(345,137){$\alpha \gamma$}
\put(335,115){\vector(1,0){40}}
\put(345,105){$\epsilon P_R$}
\put(335,15){\vector(1,0){40}}
\put(345,5){$\epsilon P_C$}
\end{picture}
\caption{The compartmental model}
\label{diagram:model}
\end{figure}
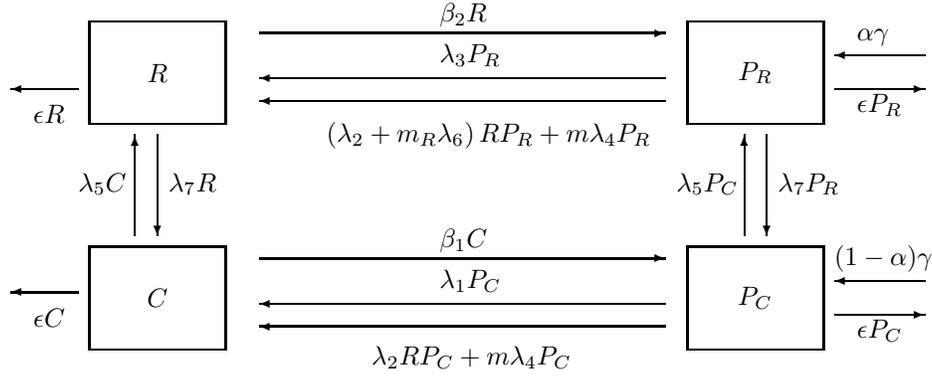

We consider a continuous compartmental model with four compartments, represented by the graph in Figure~\ref{diagram:model} and governed by an autonomous system of four ordinary differential equations.

Measuring time in years, we consider the following (pairwise disjoint) compartments: $R(t)$, the referral customers in time t, $C(t)$, the regular customers in time t, $P_R(t)$, the potential referral customers in time t and $P_C(t)$, the potential regular customers in time t. To model transitions between compartments we consider the following parameters: $\lambda_1$, the natural transition rate between $P_C$ and $C$, given by the number of potential regular customers that become regular customers without external influence per year over the number of potential customers (by ``without external influence'' we mean without being influenced by marketing campaigns or referral customers); $\lambda_2$, the referral pull effect, given by the average number of customers that a single referral brings (with no additional incentive) per year over the number of potential customers; $\lambda_3$, the natural $P_R$ to $R$ transition rate, corresponding to the number of potential referral customers that become referral customers without external influence per year over the number of potential referral customers; $m(t)$, the undifferentiated marketing costs, corresponding to marketing costs associated to undifferentiated marketing campaigns per year; $\lambda_4$, the pull effect due to undifferentiated marketing, corresponding to the quotient of the outcome of undifferentiated marketing campaigns per year by the number of potential customers (by ``outcome of undifferentiated marketing campaigns'' it is meant the number of potential customers that become customers in the sequence of undifferentiated marketing campaigns per unitary marketing cost per year); $m_R(t)$, the referral associated marketing costs, corresponding to marketing costs associated to referral directed marketing campaigns per year; $\lambda_6$, the pull effect due to referral directed marketing, given by the referral directed marketing campaigns outcome over the number of potential customers (by ``referral directed marketing campaigns outcome'' it is meant the average number of additional customers that a single referral can bring with incentives per unitary marketing cost per year); $\lambda_5$, the non-central/central transition in the social network equal to the number of individuals non-central in the social network that become central over the total number of individuals in the social network; $\lambda_7$, the central/non-central transition in the social network, given by the number of individuals central in the social network that become non-central over the total number of individuals in the social network; $\beta_1$, the regular customer defection rate, equal to the number of regular customers that cease to be customers over the number of regular customers; $\beta_2$, the referral defection rate, given by the number of referrals that cease to be customers over the number of referrals; $\epsilon$, the customer \& potential customer defection rate, corresponding to the number of individuals that leave the universe of customers and potential customers per year over the number of customers and potential customers (by ``number of individuals that leave the universe of customers and potential customers per year'' it is meant the number of customers and potential customers that cease to be in the set of customers or potential customers per year due to emigration, death, etc.); $\gamma$, the customer and potential customer recruitment rate, given by the number of individuals that enter the universe of customers and potential customers per year over the number of customers and potential customers (by ``number of individuals that leave the universe of customers and potential customers per year'' it is meant the number of customers and potential customers that cease to be in the set of customers or potential customers per year due to immigration, etc.); $\alpha$, the referral recruitment rate, equal to the number of referrals that enter the universe of customers and potential customers per year over number of individuals that enter the universe of customers and potential customers per year.

Our model can be translated into the following system of differential equations~\eqref{eq:modelo} to be studied along this paper.
{\small
\begin{equation}\label{eq:modelo}
\begin{cases}
C'= \lbd_7 R -(\eps+\beta_1+\lbd_5) C + (\lbd_1+m\lbd_4) P_C+\lbd_2RP_C\\
R'= \lbd_5 C-(\eps+\beta_2+\lbd_7)R + (\lbd_3+m\lbd_4) P_R +(\lbd_2+m_R\lbd_6)RP_R\\
P_C'=(1-\alpha)\gamma + \beta_1 C+ \lbd_7 P_R - (\eps+\lbd_5+\lbd_1+m\lbd_4) P_C - \lbd_2RP_C \\
P_R'=\alpha \gamma + \beta_2 R + \lbd_5 P_C - (\eps+\lbd_7+\lbd_3+m\lbd_4) P_R - (\lbd_2+m_R\lbd_6) RP_R
\end{cases}
\end{equation}
}
Notice that $\lambda_2 R P_R$ correspond to the average number of referrals that are brought with no incentive by the referrals per year, that $\lambda_2 R P_C$ is the average number of regular customers that are brought with no incentive by the referrals per year, that $m_R(t) \lambda_6 R P_R$ is the average number of additional customers brought due to incentives per year, that $m(t) \lambda_4 P_R$ is the number of potential referrals that become referrals in the sequence of undifferentiated marketing campaigns per year and that $m(t) \lambda_4 P_C$ is the number of potential regular customers that become regular customers in the sequence of undifferentiated marketing campaigns per year.

This paper is divided in the following way: in section~\ref{section:MR} we state our main results concerning the asymptotic behavior of the number of regular customers and referral customers, in section~\ref{section:S} we present some simulation with the objective of illustrating our theoretical results, in section~\ref{section:P} we prove our results and, finally, in section~\ref{section:Conc} we discuss the results obtained.
\section{Main Results} \label{section:MR}
One of the first natural issues to address when studying a compartmental model is the existence and stability of equilibrium solutions. We obtain several results on the existence and stability of equilibrium solutions in of model~\eqref{eq:modelo} in this section.

We first derive an auxiliary result. Given $\delta \ge 0$ define the sets
    $$\Delta_\delta=\{(x,y,z,w) \in (\R_0^+)^4: \gamma/\eps-\delta \le x+y+z+w \le \gamma/\eps+\delta\}.$$
We have the following result that shows that the total population in system~\eqref{eq:modelo} converges to the ratio $\gamma/\eps$, independently of the nonnegative initial conditions considered.

\begin{lemma}[Asymptotic behavior of the total population]\label{teo:general_system}
Let $\eps>0$ and let $(C(t),R(t),P_C(t),P_R(t))$ be some solution of system~\eqref{eq:modelo} with nonnegative initial conditions: $C(t_0) \ge 0$, $R(t_0) \ge 0$, $P(t_0) \ge 0$, $P(t_0) \ge 0$. Then:
\begin{enumerate}[a)]
  \item \label{teo:general_system-1} for all $t \ge t_0$, we have $C(t),R(t),P_C(t),P_R(t) \ge 0$;
  \item \label{teo:general_system-2} we have $\displaystyle \lim_{t \to +\infty} C(t)+R(t)+P_C(t)+P_R(t) = \dfrac{\gamma}{\eps}$.
  In particular, given $\delta>0$ and any solution $(C(t),R(t),P_C(t),P_R(t))$, with nonnegative initial conditions there is $T>0$ such that $(C(t),R(t),P_C(t),P_R(t)) \in \Delta_\delta$ for all $t \ge T$ and any equilibrium solution is in the set
  $\Delta_0=\{(x,y,z,w) \in (\R_0^+)^4: x+y+z+w = \gamma/\eps\}$.
\end{enumerate}

\end{lemma}

\begin{remark}
  The case $\eps=0$ is not a very interesting case since it corresponds to the situation where there is no customer \& potential customer defection rate which is not a realistic assumption. Nevertheless it is easy to check that, if $\eps=0$, then, given initial conditions $C_0$, $R_0$, $P_{C,0}$ and $P_{R,0}$, we have
    $$\displaystyle \lim_{t \to +\infty} C(t)+R(t)+P_C(t)+P_R(t) = \gamma t  + C_0+R_0+P_{C,0}+P_{R,0}.$$
  In particular if $\eps=\gamma=0$, then the total population remains constant.
\end{remark}

 We now obtain a result on the existence of equilibrium solutions. Under  the assumption of positivity of the defection rate, the referral pull effect and the non central/central transition in the social network we conclude that there are one, two or three equilibrium solutions, depending on the number of real roots of some third degree polynomial. We need to define the constants
\begin{equation}\label{eq:p_q}
  p=\dfrac{\gamma(\alpha\eps+\lbd_5)}{\eps(\eps+\lbd_5+\lbd_7)} \quad \quad \text{and}
  \quad \quad q=\dfrac{\gamma((1-\alpha)\eps+\lbd_7)}{\eps(\eps+\lbd_5+\lbd_7)}
\end{equation}
and also
\[
  u=\eps+\beta_1+\lambda_5+\lambda_1+m\lambda_4
  \quad \text{and} \quad
  v=\eps+\beta_2+\lambda_7+\lambda_3+m\lambda_4.
\]

\begin{theorem}[Equilibrium solutions]\label{teo:equilibriums}
Let $\eps,\lambda_5>0$. Then:
\begin{enumerate}[a)]
\item \label{teo:equilibriums-1} system~\eqref{eq:modelo} has up to three equilibrium solutions $(R^*,C^*,P_R^*,P_C^*)$. The first component, $R^*$, is always a nonnegative solution of the cubic equation
\begin{equation}\label{eq:solution-R}
aR^3+bR^2+cR+d=0,
\end{equation}
where $a=-\lambda_2(\lambda_2+m_R\lambda_6)$, $b=\lbd_2(\lbd_2+m_R\lbd_6)p-u(\lambda_2+m_R\lambda_6)-\lambda_2 v$, $c=\lbd_2(\lbd_3+m\lbd_4)p+u(\lbd_2+m_R\lbd_6)p+(\lbd_7+\lbd_2 q)\lbd_5-uv$ and $d=(\lambda_3+m\lambda_4)pu+(\lambda_1+m\lambda_4)q\lambda_5$;
\item \label{teo:equilibriums-2} any equilibrium solution $(R^*,C^*,P_R^*,P_C^*)$ is obtained in the following way: $R^*$ is a nonnegative solution of~\eqref{eq:solution-R} and $P_R^*=p-R^*$, $P_C^*=q-C^*$ and
  \begin{equation}\label{eq:C*}
  C^*=\dfrac{(\eps+\beta_2+\lambda_7)R^*- (\lambda_3+m\lambda_4+(\lambda_2+m_R\lambda_6)R^*)(p-R^*)}{\lambda_5},
  \end{equation}
   are nonnegative constants.
\end{enumerate}
\end{theorem}

In the next result we discuss the asymptotic behavior of solutions of~\eqref{eq:modelo} under some assumptions on the parameters that roughly correspond to require that the referral pull effects are bounded by some functions that we can identify with the other ``forces'' in the model such as the natural transition rates, the pull effects due to undifferentiated marketing and the defection rates (see equation~\eqref{eq:cond-equiv}).
In the following theorem we were able to show that, under the mentioned assumptions, the asymptotic behavior of the solutions of~\eqref{eq:modelo} can be obtained by the two dimensional autonomous system~\eqref{eq:sistema_equiv_para_C_and_R_=_0}.

\begin{theorem}[Asymptotic behavior of solutions] \label{teo:asymptotic_behavior}
Let $\eps>0$ and assume that
\begin{equation}\label{eq:cond-equiv}
\small \min\left\{\dfrac{2\eps+2\beta_1+\lbd_5+2\lbd_1+2m\lbd_4}{\lbd_7+q\lbd_2}, \, \frac{2\eps+2\beta_2+\lbd_7+2\lbd_3+2m\lbd_4+2(\lbd_2+m_R\lbd_6)p}{\lbd_5+q\lbd_2} \right\}\!>\!1.
\end{equation}
Consider the system
          {\small
          \begin{equation}\label{eq:sistema_equiv_para_C_and_R_=_0}
            \begin{cases}
                C_a'=\lbd_7 R_a -(\eps+\beta_1+\lbd_5) C_a + (\lbd_1+m\lbd_4+\lbd_2 R_a) (q-C_a)\\
                R_a'=\lbd_5 C_a-(\eps+\beta_2+\lbd_7)R_a + (\lbd_3+m\lbd_4+(\lbd_2+m_R\lbd_6)R_a) \left(p - R_a\right)
            \end{cases}
          \end{equation}}
and set $P_{R,a}(t)=p-R_a(t)$ and $P_{C,a}(t)=q-C_a(t)$. Then the asymptotic behavior of $C$, $R$, $P_C$ and $P_R$ in system~\eqref{eq:modelo} is the same as the asymptotic behavior of $C_a$, $R_a$, $P_{C,a}$ and $P_{R,a}$. Namely if $(C(t),R(t),P_C(t),P_R(t))$ is a solution of~\eqref{eq:modelo} with initial condition $(C(t_0),R(t_0),P_C(t_0),P_R(t_0))=(C_0,R_0,P_{C,0},P_{R,0})$ and  $(C_a(t),R_a(t))$ is a solution of~\eqref{eq:sistema_equiv_para_C_and_R_=_0} with initial condition $(C_a(t_0),R_a(t_0))=(C_0,R_0)$ then
$$\lim_{t \to +\infty} |C(t)-C_a(t)|+|R(t)-R_a(t)|+|P_C(t)-P_{C,a}(t)|+|P_R(t)-P_{R,a}(t)|=0.$$
\end{theorem}

In the next two theorems we discuss two particular situations where we analyse the existence of equilibriums and their stability. The vector fields plotted with the objective of illustrating the situations correspond to the reduced system~\eqref{eq:sistema_equiv_para_C_and_R_=_0} but we considered situations where Theorem~\ref{teo:asymptotic_behavior} apply so that the asymptotic behavior of referrals and regular clients is the same for both systems.

First we will discuss the situation where there is no transition between referral/potential referral and customer/potential customer and thus we set $\lbd_5=\lbd_7=0$. We consider two cases: the situation where $\lbd_3+m\lbd_4 \ne 0$ and $\lbd_1+m\lbd_4 \ne 0$ (we named it ``static social network'' to reflect the fact that there is no transition between referrals and regular customers) and the situation where $\lbd_1=\lbd_2=\lbd_4=0$, corresponding to the case where all potential customers and potential referrals that become customers are consequence of referral influence (we named it ``word of mouth'' to emphasise that all marketing efforts are related to referrals).

We have the following result in the static social network case.

\begin{theorem}[Static social network]\label{teo:static_social_network}
 The following holds for system~\eqref{eq:modelo} with $\eps>0$, $\lbd_2>0$, $\lbd_3+m\lbd_4 \ne 0$, $\lbd_1+m\lbd_4 \ne 0$ and $\lbd_5=\lbd_7=0$: there is a unique equilibrium solution $(R^*,C^*,P_R^*,P_C^*)$ that is locally asymptotically stable and is given by
    \begin{equation}\label{eq:static-R*}
    R^* = \dfrac{\alpha\gamma}{2\eps}-\theta+\sqrt{\left(\dfrac{\alpha\gamma}{2\eps}+\theta\right)^2 -\frac{\alpha\gamma(\eps+\beta_2)}{\eps(\lbd_2+m_R\lbd_6)}},
    \end{equation}
    $$C^* = \dfrac{(1-\alpha)\gamma(\lbd_1+m\lbd_4+\lbd_2 R^*)}{\eps(\eps+\beta_1+\lbd_1+m\lbd_4+\lbd_2 R^*)},$$\\[1mm]
    $P_R^*=\gamma\alpha/\eps-R^*$ and $P_C^*=\gamma(1-\alpha)/\eps-C^*$ where
\[
  \theta=\dfrac{\lbd_3+m\lbd_4+\beta_2+\eps}{2(\lbd_2+m_R \lbd_6)}.
\]
\end{theorem}

Define
\begin{equation}\label{eq:tau}
  \tau=\dfrac{\alpha\gamma(\lbd_2+m_R\lbd_6)}{\eps(\eps+\beta_2)}.
\end{equation}

\begin{figure}[h!]
    \begin{minipage}[b][4cm]{.5\linewidth}
    \includegraphics[scale=0.37]{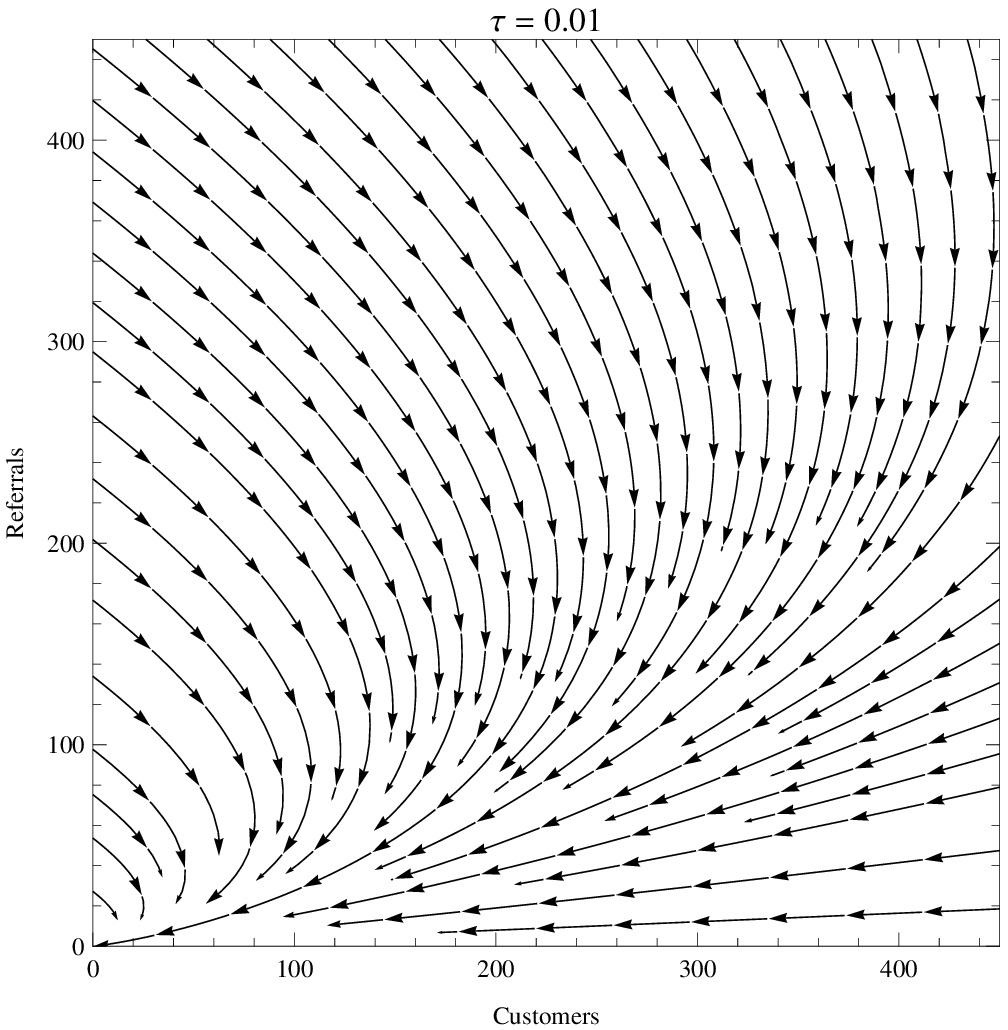}
  \end{minipage}
  \begin{minipage}[b][4cm]{.2\linewidth}
    \includegraphics[scale=0.37]{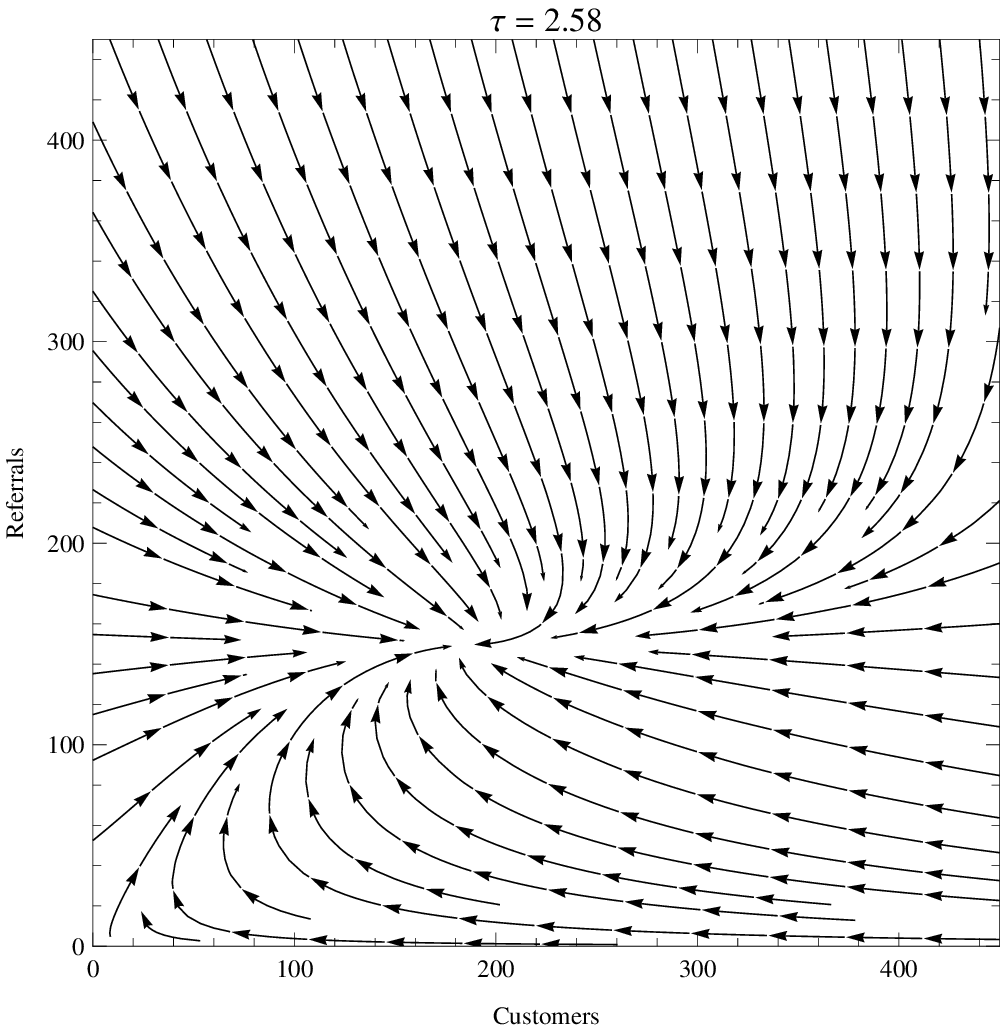}
  \end{minipage}
    \caption{Word of mouth: $\tau<1$ and $\tau>1$}
      \label{fig1}
\end{figure}

We now consider the word of mouth case. In Figure~\ref{fig1} we show the behavior in the plane $C-R$ in both of the regimes of Theorem~\ref{teo:mri}.
In Figure~\ref{fig1}  we used for the plot in the left $m = 40$, $m_R = 0$, $\eps=0.01$ and $\lambda_2 = 10^{-5}$, and, for the plot in the right, $m = 30$, $m_R = 10$, $\eps=0.01$ and $\lambda_2 = 10^{-5}$.

\begin{theorem}[Word of mouth]\label{teo:mri}
 The following statements holds for system~\eqref{eq:modelo} with $\eps>0$, $\lbd_2>0$ and $\lbd_1=\lbd_3=\lbd_4=\lbd_5=\lbd_7=0$:
 \begin{enumerate}[a)]
   \item if $\tau \le 1$ then there a unique locally stable equilibrium given by
       \begin{equation*}
       (C^*,R^*,P_C^*,P_R^*) =\left(0,0,\frac{\gamma(1-\alpha)}{\eps},\frac{\gamma\alpha}{\eps}\right);
       \end{equation*}
   \item if $\tau > 1$ then there are two equilibrium solutions. An unstable equilibrium given by
       \begin{equation*}
       (C_1^*,R_1^*,P_{C,1}^*,P_{R,1}^*) =\left(0,0,\frac{\gamma(1-\alpha)}{\eps},\frac{\gamma\alpha}{\eps}\right);
       \end{equation*}
       and a locally stable equilibrium given by
       \begin{equation*}
       \begin{split}
       & (C_2^*, R_2^*, P_{C,2}^*,P_{R,2}^*) \\
       & = \left(\frac{\alpha(1-\alpha)\lambda_2\gamma^2(1-\frac{1}{\tau})} {\eps^2(\eps+\beta_1+\lbd_2(1-\frac{1}{\tau}))},
       \dfrac{\alpha\gamma(1-\frac{1}{\tau})}{\eps},
       \frac{(1-\alpha)(\eps+\beta_1)\gamma}{\eps(\eps+\beta_1+\lbd_2(1-\frac{1}{\tau}))}, \frac{\eps+\beta_2}{\lbd_2+m_R\lbd_6}\right).
       \end{split}
       \end{equation*}
 \end{enumerate}
\end{theorem}

We next consider a scenario where there is no direct referral influence and thus we set $\lbd_2=\lbd_6=0$. Define
\[
  \kappa_1=\dfrac{\lbd_5[\lbd_7p+(\lbd_1+m\lbd_4)q]}{up(\eps+\beta_2+\lbd_7)}
  \quad \text{and} \quad \kappa_2=\dfrac{\lbd_7[\lbd_5q+(\lbd_3+m\lbd_4)p]}{vq(\eps+\beta_1+\lbd_5)}.
\]

\begin{theorem}[No referral influence]\label{teo:no_referral_influence}
 For system~\eqref{eq:modelo} with $\eps>0$, $\lbd_5>0$ and $\lbd_2=\lbd_6=0$ there is a unique equilibrium solution that is globally asymptotically stable and is given by
        $$(R^*, C^*, P_R^*, P_C^*)=(\kappa_1 p, \, \kappa_2 p, \, p-\kappa_1p, \, q-\kappa_2p).$$
\end{theorem}

\begin{figure}[h!]
    \begin{minipage}[b][4cm]{.5\linewidth}
    \includegraphics[scale=0.37]{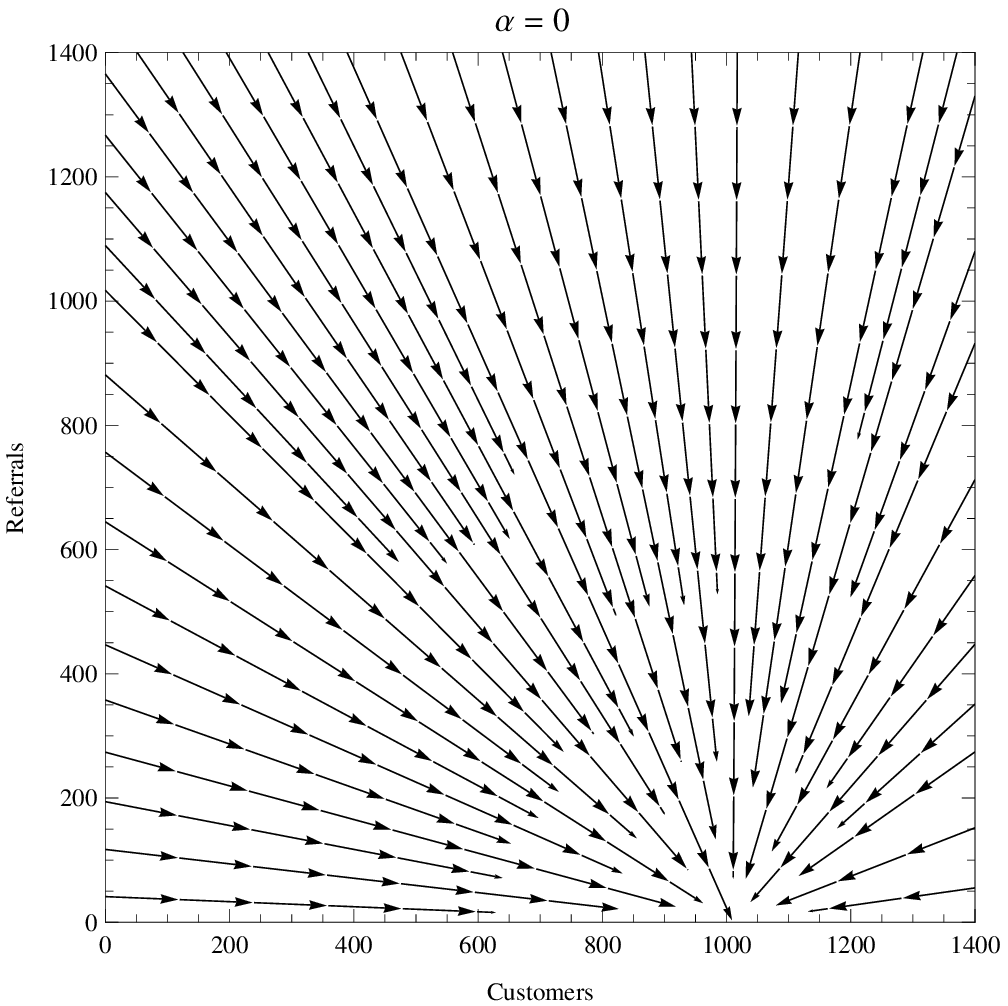}
  \end{minipage}
    \begin{minipage}[b][4cm]{.2\linewidth}
    \includegraphics[scale=0.37]{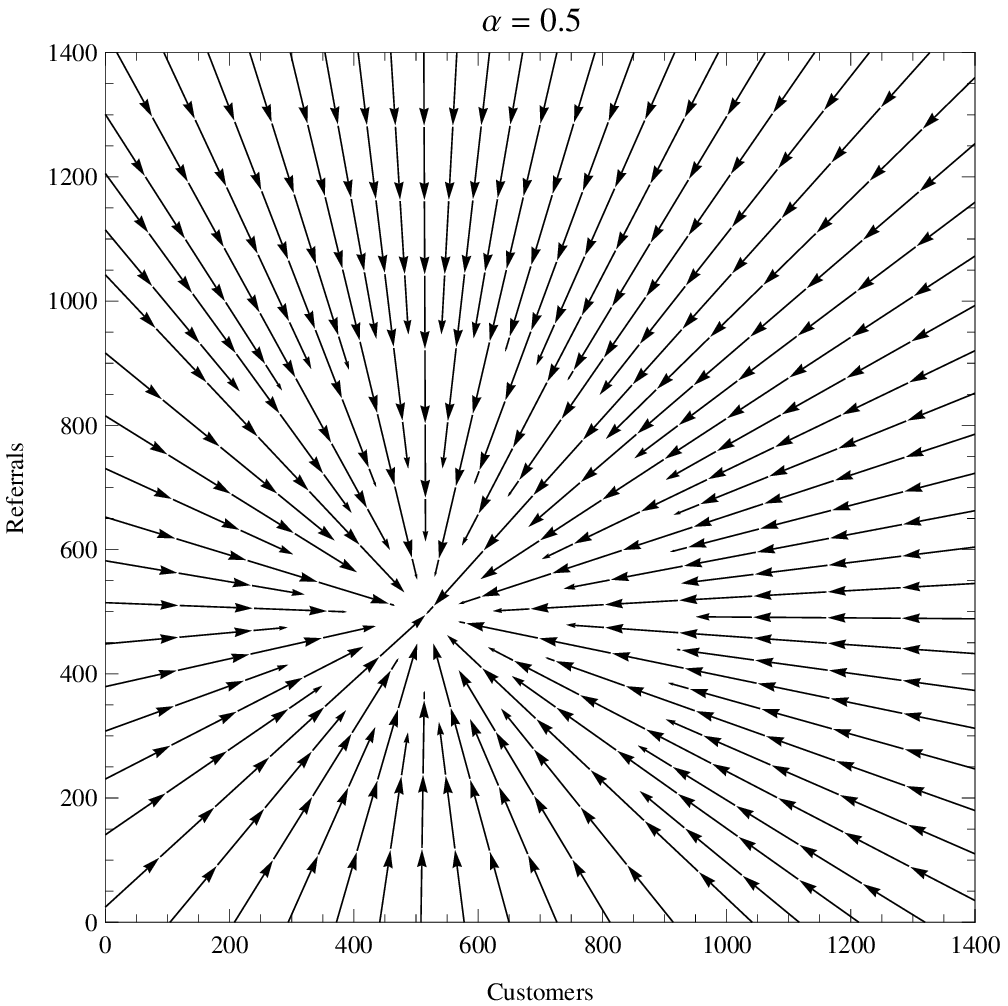}
  \end{minipage}
    \caption{Word of mouth: $\alpha=0$ and $\alpha=0.5$}
      \label{fig2}
\end{figure}

In Figure~\ref{fig2} we illustrate the behavior of the reduced system in the plane $C-R$ for the setting in Theorem~\ref{teo:no_referral_influence}. In Figure~\ref{fig2} we used $\lbd_2=\lbd_6=0$, $m = 40$, $m_R = 0$, $\lbd_1=\lbd_3=\lbd_4=\lbd_7=0.0002$, $\lbd_5=0.0000018$, $\beta_1=\beta_2=0.18$, $\eps=0.01$ and, for the plot on the left, $\alpha=0$, and, for the plot on the right $\alpha=0.5$. Note that for the figure in the left the number of referrals in the equilibrium point is nonzero, although it seems to be (the number of referrals is very low since the there are no referrals entering the population).

\section{Simulation} \label{section:S}

To obtain a better understanding of the behavior of our model, we assume that, for a given corporation, we consider the values for the parameters presented in Table~\ref{tab:ref}.
\begin{table}[htb]
\centering
{\renewcommand{\arraystretch}{1.5} 
\begin{tabular}{cc}
\toprule
Parameter & Value \\
\midrule
$\alpha$, $\epsilon$& 0.01\\
$\gamma$ & $\epsilon N_0$\\
$\lambda_1,\ldots,\lambda_4,\lambda_6, \lambda_7$ & 0.0002 \\
$\lambda_5$ & $\lambda_7 R_0/C_0$\\
$\beta_1,~ \beta_2$ & 0.18\\
\bottomrule
\end{tabular}}
\caption{Values of parameters}
\label{tab:ref}
\end{table}
and initial conditions $C_0=2200$, $P_{C,0}=22000$, $R_0=20$ and $P_{R,0}=200$. We write $N_0=C_0+P_{C,0}+R_0+P_{R,0}$. The value considered for $\beta_1$ and $\beta_2$ is based on usual assumptions concerning the defection rate~\cite{Lee-Lee-Feick-JDBCSM-2006}. The assumptions $\lbd_5 C_0=\lbd_7 R_0$ and $\lbd_5 P_{C,0}=\lbd_7 P_{R,0}$ are made to assure that the underlying social network hasn't an initial tendency to ``benefit'' any of the four compartments. We also consider $\gamma=N_0\times \epsilon$ so that the total population converges to an equilibrium where the total population equals $N_0$.

We solved system~\eqref{eq:modelo} (named \emph{initial}) and system~\eqref{eq:sistema_equiv_para_C_and_R_=_0} (named \emph{reduced}) with MATLAB. In the figures we plot the solution for the initial and the reduced systems.  We considered two sets of values for $m$ and $m_R$, namely $(m,m_R)=(40,0)$ corresponding to a situation of undifferentiated marketing and $(m,m_R)=(30,10)$ corresponding to a situation where some effort is made for attracting referrals. In both situations we maintain the same total effort, $m+m_R=40$, in order to be able to compare both cases.

\begin{figure}[!htb]
        \centering
        \begin{subfigure}[b]{0.49\textwidth}
                \includegraphics[scale=0.4]{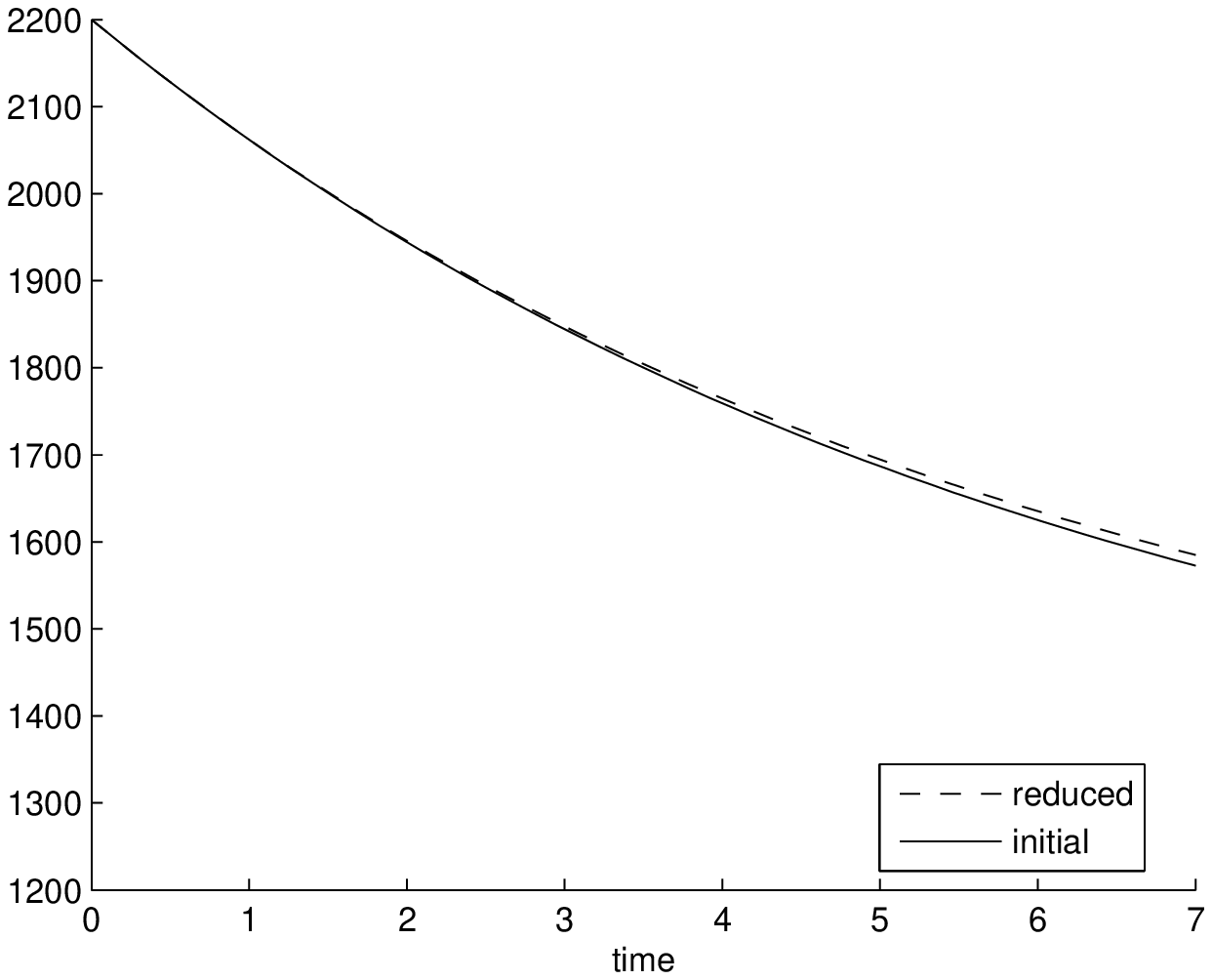}
                \caption{$t\in[0;7]$}
        \end{subfigure}%
        ~ 
        \begin{subfigure}[b]{0.49\textwidth}
                \includegraphics[scale=0.4]{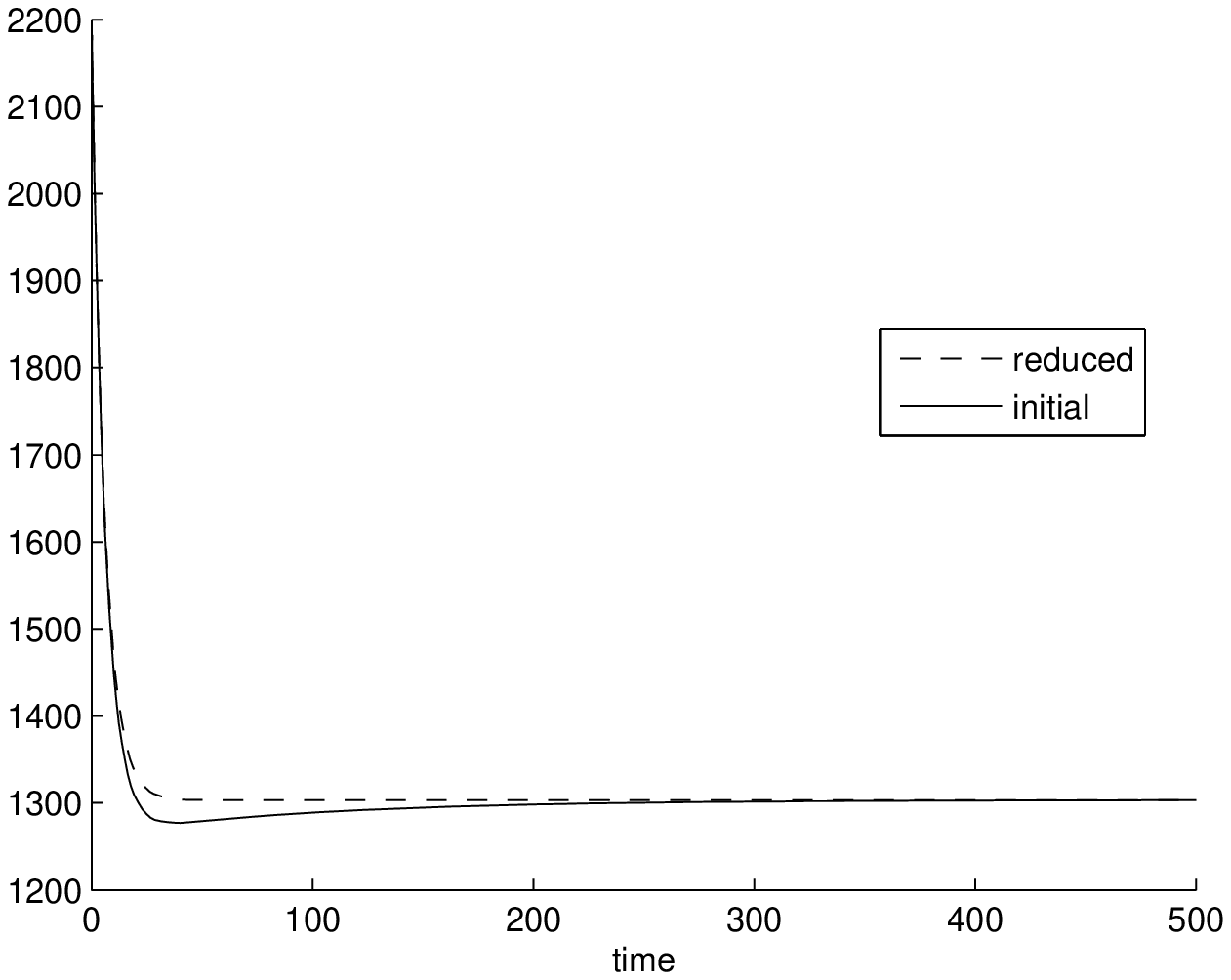}
                \caption{$t\in[0;500]$}
        \end{subfigure}
        \caption{Regular customers evolution with $m=40$ and $m_R=0$}
        \label{fig:undif-costumers}
\end{figure}

\begin{figure}[!htb]
        \centering
        \begin{subfigure}[b]{0.49\textwidth}
                \includegraphics[scale=0.4]{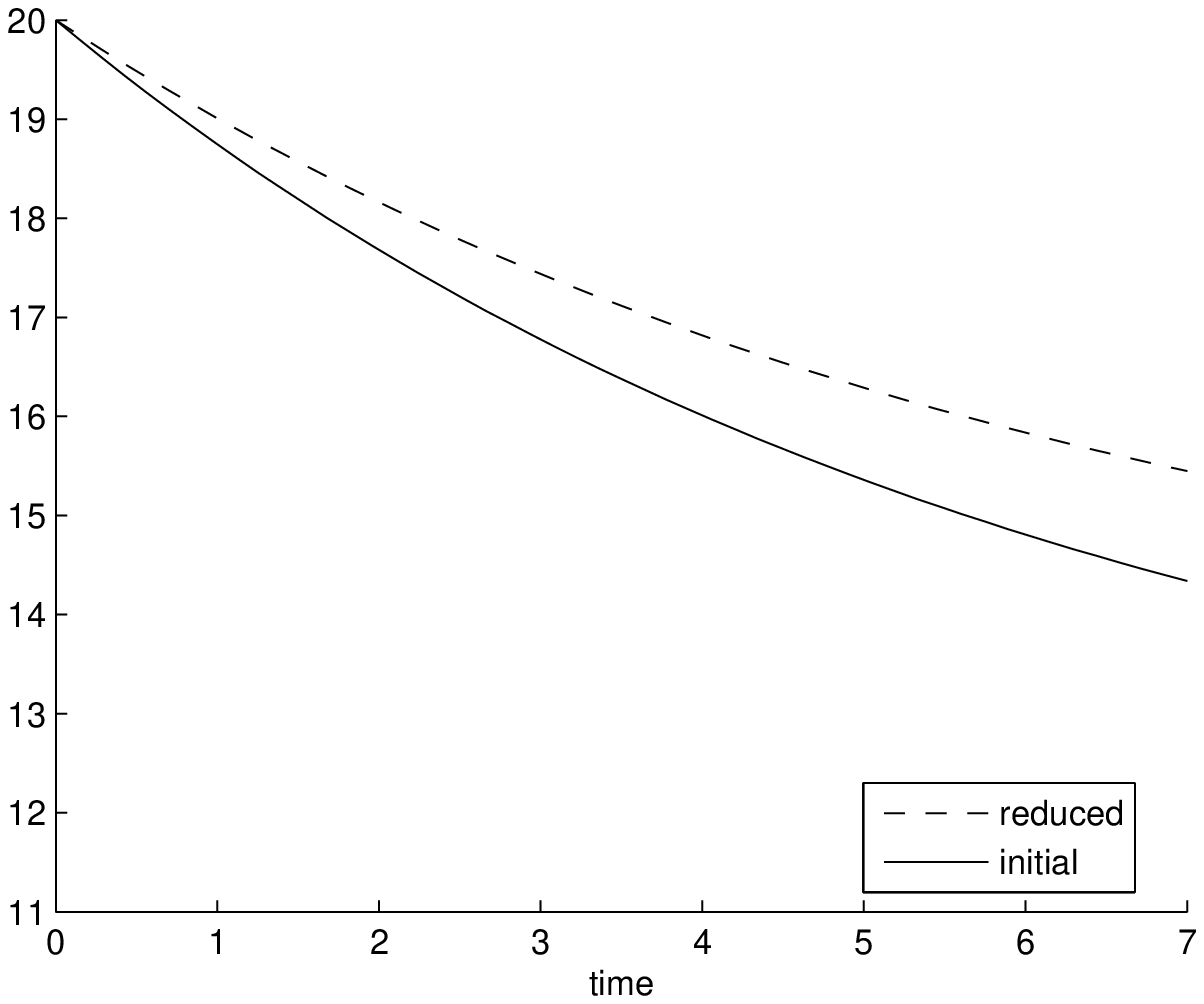}
                \caption{$t\in[0;7]$}
        \end{subfigure}%
        ~ 
        \begin{subfigure}[b]{0.49\textwidth}
                \includegraphics[scale=0.4]{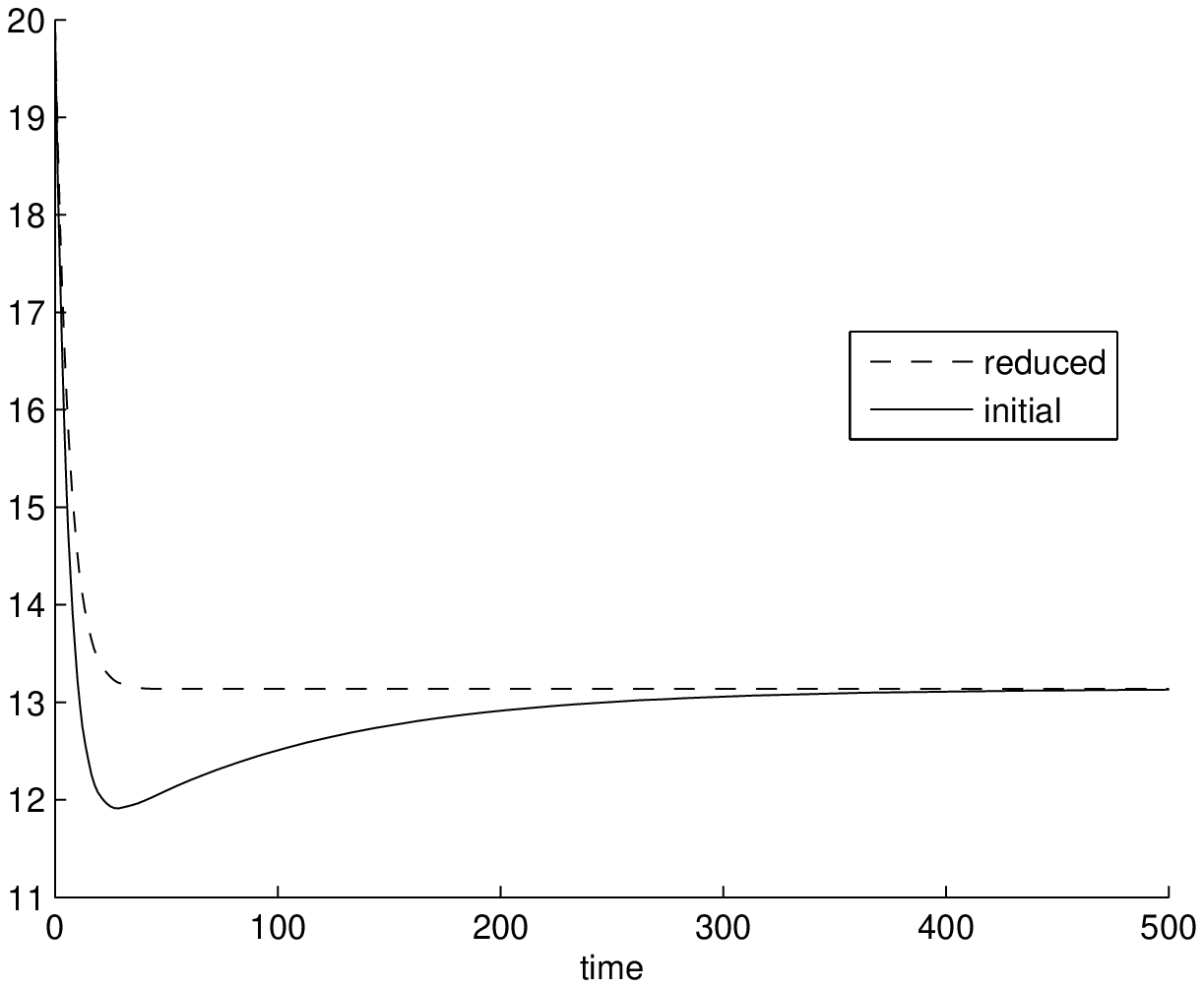}
                \caption{$t\in[0;500]$}
        \end{subfigure}
        \caption{Referral evolution with $m=40$ and $m_R=0$}
        \label{fig:undif-referrals}
\end{figure}

In Figures~\ref{fig:undif-costumers} and~\ref{fig:undif-referrals} we consider the evolution of customers and referrals in the case where marketing is used in an undifferentiated way. We can see that the number of customers and referrals decreases in this situation and stabilizes in some lower value for both compartments.

\begin{figure}[!htb]
        \centering
        \begin{subfigure}[b]{0.49\textwidth}
                \includegraphics[scale=0.4]{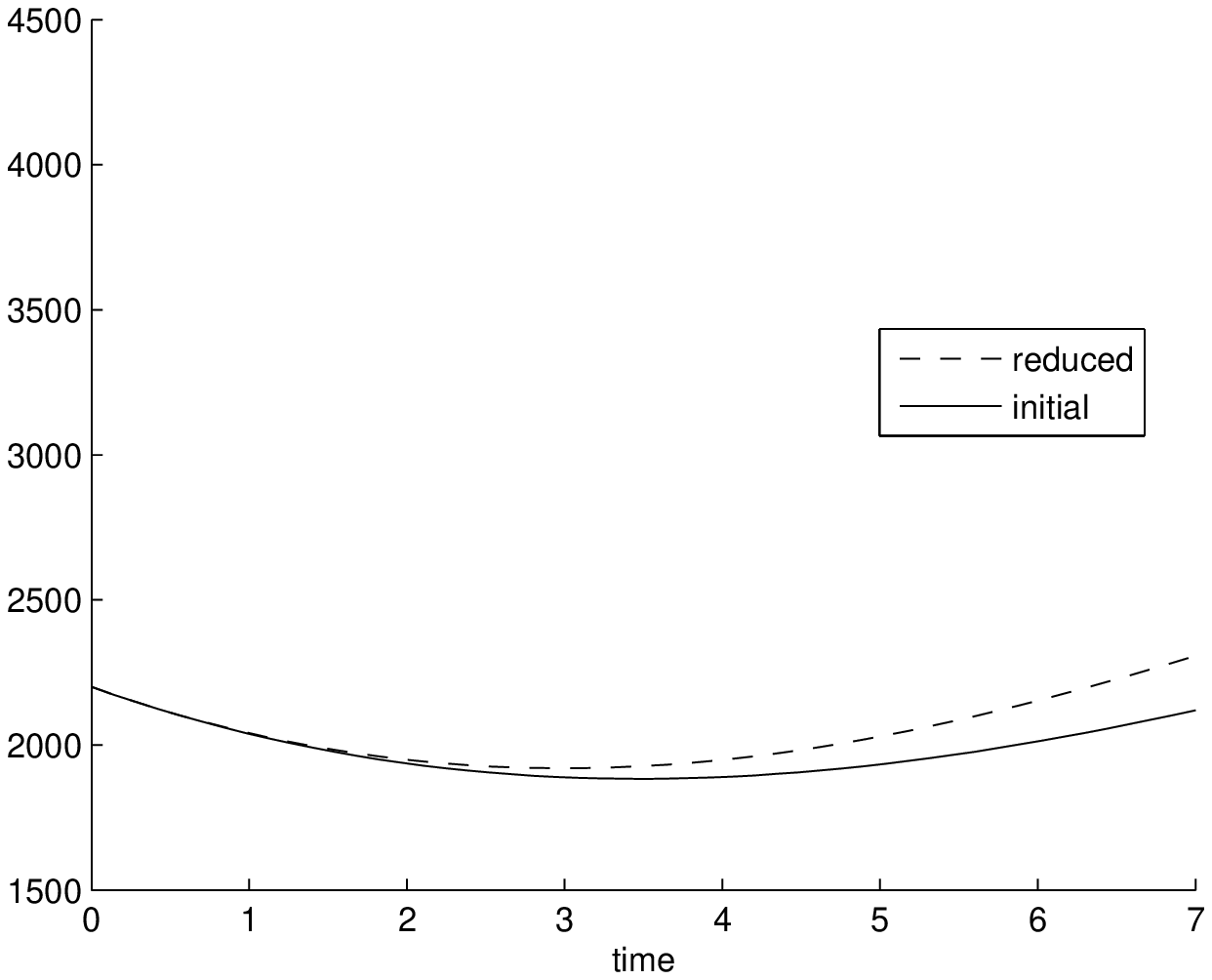}
                \caption{$t\in[0;7]$}
        \end{subfigure}%
        ~ 
        \begin{subfigure}[b]{0.49\textwidth}
                \includegraphics[scale=0.4]{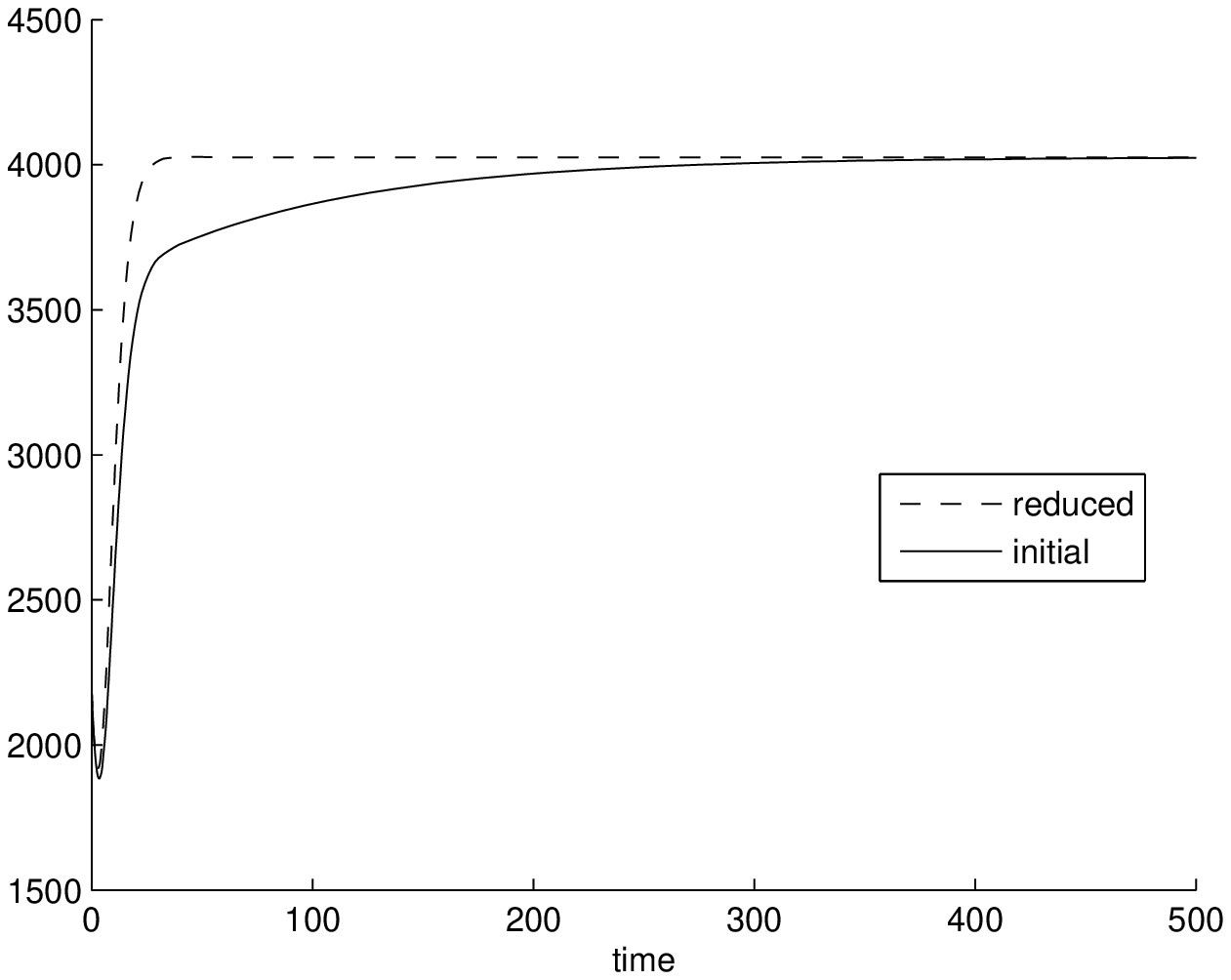}
                \caption{$t\in[0;500]$}
        \end{subfigure}
        \caption{Regular customers evolution with $m=30$ and $m_R=10$}
        \label{fig:mark-costumers}
\end{figure}

\begin{figure}[!htb]
        \centering
        \begin{subfigure}[b]{0.49\textwidth}
                \includegraphics[scale=0.4]{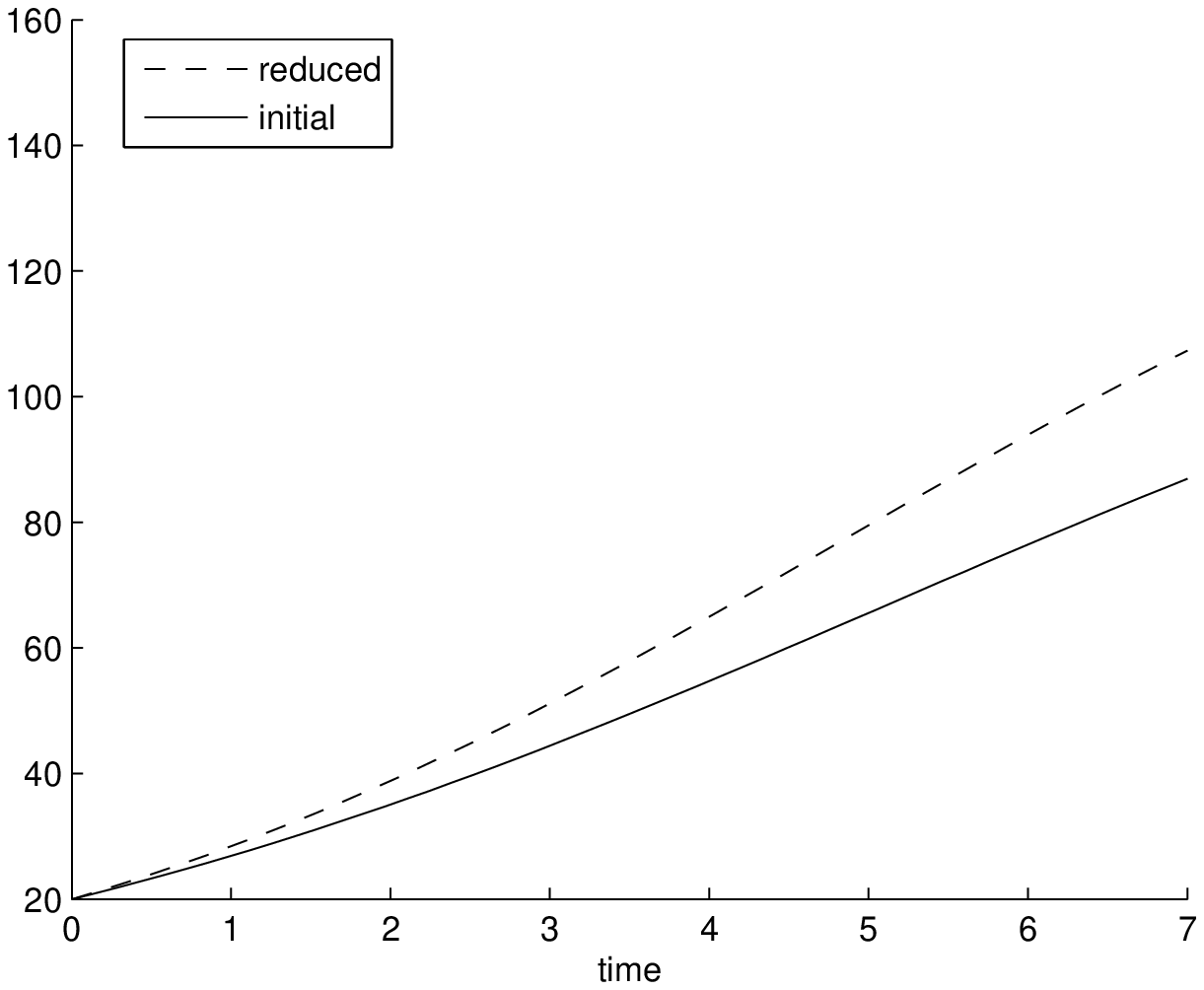}
                \caption{$t\in[0;7]$}
        \end{subfigure}%
        ~ 
        \begin{subfigure}[b]{0.49\textwidth}
                \includegraphics[scale=0.4]{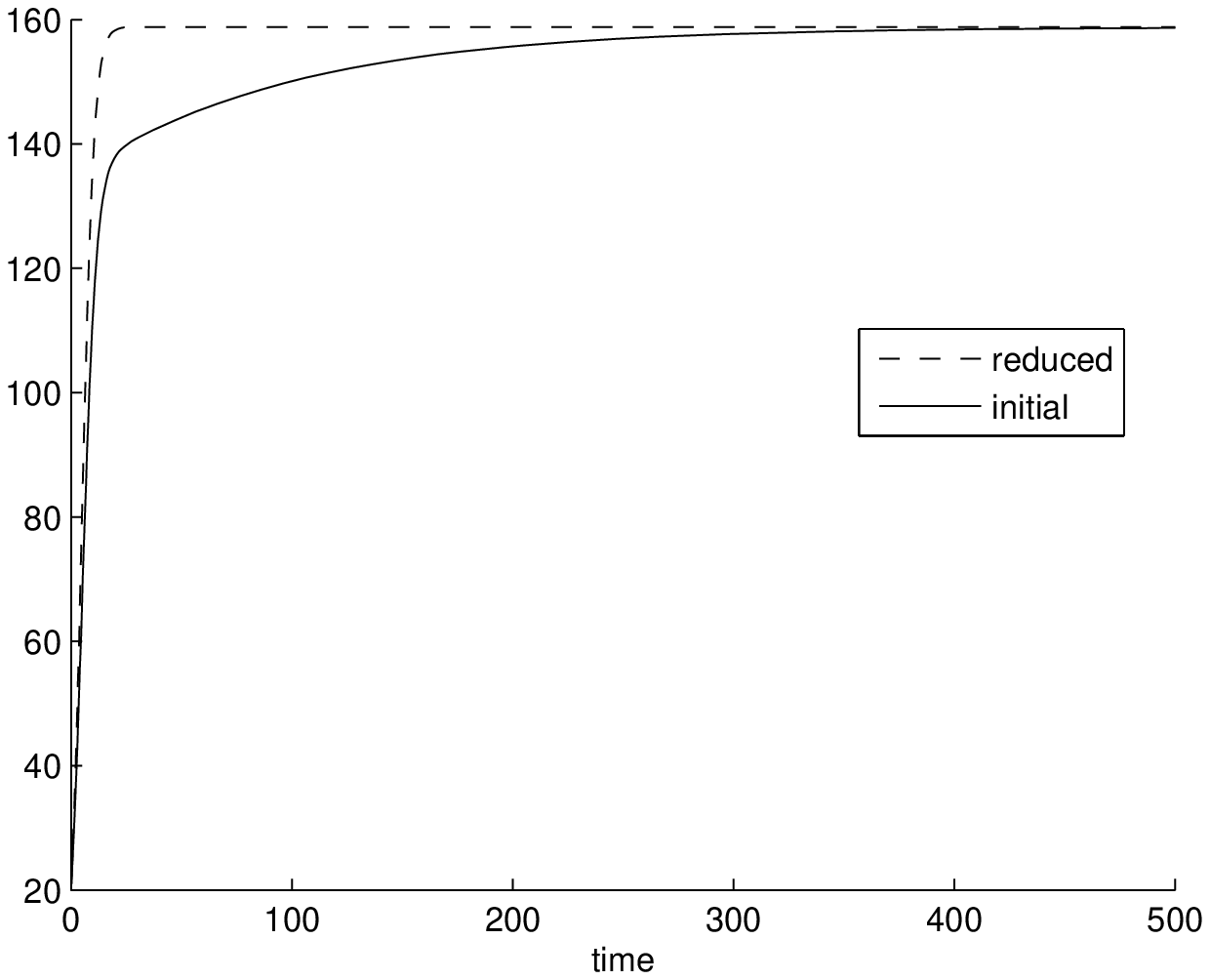}
                \caption{$t\in[0;500]$}
        \end{subfigure}
        \caption{Referral evolution with $m=30$ and $m_R=10$}
        \label{fig:mark-referrals}
\end{figure}

In Figures~\ref{fig:mark-costumers} and~\ref{fig:mark-referrals} we now consider the evolution of customers and referrals in the case where some marketing effort is used to attract referrals. We can see that there is an initial small decrease in the number of customers that is followed rapidly by an increase that asymptotically doubles its number. There is also an increase in the number of referrals due to the positive value of $m_R$.

\begin{figure}[!htb]
        \centering
        \begin{subfigure}[b]{0.49\textwidth}
                \includegraphics[scale=0.4]{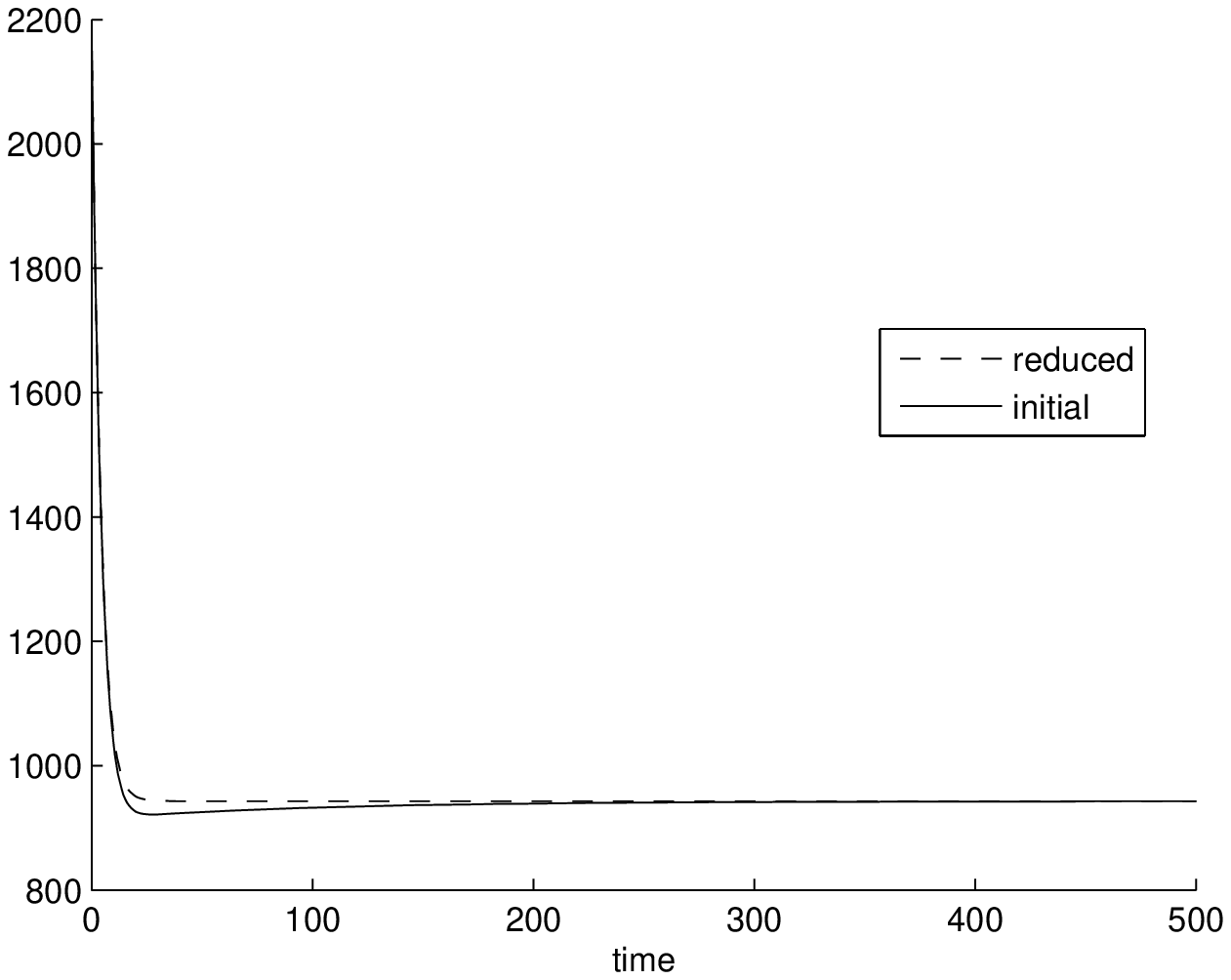}
                \caption{Regular customers evolution}
        \end{subfigure}%
        ~ 
        \begin{subfigure}[b]{0.49\textwidth}
                \includegraphics[scale=0.4]{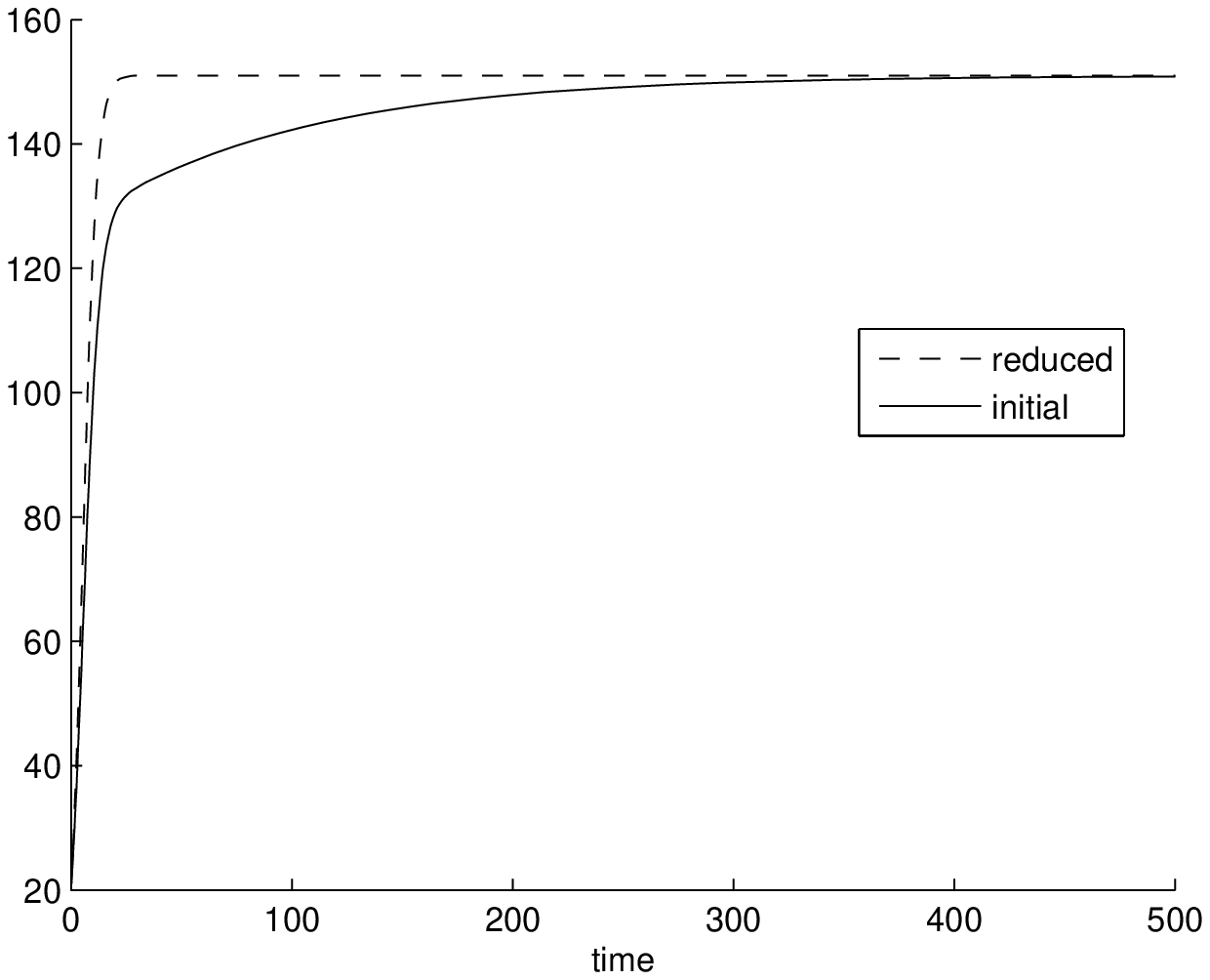}
                \caption{Referrals evolution}
        \end{subfigure}
        \caption{Asymptotic behavior with $m=30$,  $m_R=10$ and $\lambda_2=1\times 10^{-5}$}
        \label{fig:asymptotic}
\end{figure}

In the Figure~\ref{fig:asymptotic} we present the evolution of customers and referrals in the case where $\lbd_2$ is reduced to $10^{-5}$ in order to satisfy the condition~\eqref{eq:cond-equiv} in Theorem~\eqref{teo:asymptotic_behavior}. We can see that, as stated in the theorem, the solutions asymptotically approach the same value. In the previous cases, although condition~\eqref{eq:cond-equiv} is not satisfied, there is computational evidence that the same happens. Thus we conjecture that condition~\eqref{eq:cond-equiv} can be weakened.

\vspace{1cm}

\section{Proofs} \label{section:P}

\subsection{Proof of Lemma~\ref{teo:general_system}}
Let $\eps >0$. Analysing the direction of the flow on the boundary of the set $\{(x,y,z,w) \in \R^4: x,y,z,w \ge 0\}$ we immediately obtain~\ref{teo:general_system-1}). Adding the four equations in~\eqref{eq:modelo} and letting $N(t)=C(t)+R(t)+P_C(t)+P_R(t)$, we get
the differential equation
\begin{equation}\label{eq:total pop}
N'=\gamma-\eps N.
\end{equation}
The general solution of~\eqref{eq:total pop}
is $N(t)=\gamma/\eps + C \e^{-\eps t}$ and thus $\displaystyle \lim_{t \to +\infty} N(t) = \gamma/\eps$. In particular, if $N(t)$ corresponds to an equilibrium solution then $N(t)=\gamma/\eps$. This proves~\ref{teo:general_system-2}).
\subsection{Proof of Theorem~\ref{teo:equilibriums}}
Adding equations for $R$ and $P_R$ and for $C$ and $P_C$ in~\eqref{eq:modelo} we get
\[
\begin{cases}
(R+P_R)'=\alpha\gamma-(\eps+\lbd_7)(R+P_R) + \lbd_5 (C+P_C)\\
(C+P_C)'=(1-\alpha)\gamma-(\eps+\lbd_5) (C+P_C) + \lbd_7 (R+ P_R)
\end{cases}
\]
and, setting $X=R+P_R$ and $Y=C+P_C$, we obtain
\begin{equation}\label{eq:sist_X_Y}
\begin{cases}
X'=\alpha\gamma-(\eps+\lbd_7) X + \lbd_5 Y \\
Y'=(1-\alpha)\gamma+\lbd_7 X -(\eps+\lbd_5) Y
\end{cases}.
\end{equation}
The linear system~\eqref{eq:sist_X_Y} has the general solution
\[
\begin{cases}
X(t)=\dfrac{\gamma(\alpha\eps+\lbd_5)}{\eps(\eps+\lbd_5+\lbd_7)}+\dfrac{C_1 \lbd_7 - C_2 \lbd_5}{\lbd_5+\lbd_7}\e^{-(\eps+\lbd_5+\lbd_7)t}+\dfrac{(C_1+C_2)\lbd_5}{\lbd_5+\lbd_7} \e^{-\eps t} \\
Y(t)=\dfrac{\gamma((1-\alpha)\eps+\lbd_7)}{\eps(\eps+\lbd_5+\lbd_7)}+\dfrac{C_2 \lbd_5 - C_1 \lbd_7}{\lbd_5+\lbd_7}\e^{-(\eps+\lbd_5+\lbd_7)t}+\dfrac{(C_1+C_2)\lbd_7}{\lbd_5+\lbd_7} \e^{-\eps t}
\end{cases}.
\]
and thus
\begin{equation}\label{eq:elim_Pr_Pc}
\begin{cases}
P_R(t)=\dfrac{\gamma(\alpha\eps+\lbd_5)}{\eps(\eps+\lbd_5+\lbd_7)}+\theta_1(t) - R(t)\\
P_C(t)=\dfrac{\gamma((1-\alpha)\eps+\lbd_7)}{\eps(\eps+\lbd_5+\lbd_7)}+\theta_2(t) - C(t)
\end{cases}
\end{equation}
with $\theta_i(t) \to 0$ as $t \to \infty$ for $i=1,2$. Replacing $P_R$ and $P_C$ by the expressions in~\eqref{eq:elim_Pr_Pc} in the first two equations in system~\eqref{eq:modelo} we obtain
\begin{equation}\label{eq:C_and_R}
\begin{cases}
C'=\lbd_7 R -(\eps+\beta_1+\lbd_5) C + (\lbd_1+m\lbd_4+\lbd_2R) \left( q+\theta_2(t) - C \right)\\
R'=\lbd_5 C-(\eps+\beta_2+\lbd_7)R + (\lbd_3+m\lbd_4+(\lbd_2+m_R\lbd_6)R) \left(p+\theta_1(t) - R\right)
\end{cases},
\end{equation}
where $p$ and $q$ are given by~\eqref{eq:p_q}.

By Lemma~\ref{teo:general_system} every equilibrium solution must belong to $\Delta_0$. Thus every equilibrium solution must satisfy
\begin{equation}\label{eq:C_and_R_equilibrium}
\begin{cases}
\lbd_7 R -(\eps+\beta_1+\lbd_5) C + (\lbd_1+m\lbd_4+\lbd_2R) \left( q - C \right)=0\\
\lbd_5 C-(\eps+\beta_2+\lbd_7)R + (\lbd_3+m\lbd_4+(\lbd_2+m_R\lbd_6)R) \left(p - R\right)=0
\end{cases}.
\end{equation}
Since $\lambda_5 \ne 0$ we have
\begin{equation}\label{eq:C}
C=\dfrac{(\eps+\beta_2+\lbd_7)R - (\lbd_3+m\lbd_4+(\lbd_2+m_R\lbd_6)R) \left(p - R\right)}{\lbd_5}.
\end{equation}
Substituting~\eqref{eq:C} in the first equation of~\eqref{eq:C_and_R_equilibrium} we obtain equation~\eqref{eq:solution-R}. Each solution of~\eqref{eq:solution-R} uniquely determines $C$, $R$ and $P_R$. Thus, we have at most three equilibrium solutions and~\ref{teo:equilibriums-1}) follows.

By~\eqref{eq:elim_Pr_Pc} for any equilibrium solution $(R^*,C^*,P_R^*,P_C^*)$ we have $P_R^*=p-R^*$ and $P_C^*=q-C^*$. By~\eqref{eq:C} we obtain~\eqref{eq:C*} and~\ref{teo:equilibriums-2}).
\subsection{Proof of Theorem~\ref{teo:asymptotic_behavior}}
Assume that $\eps>0$ and choose $\delta>0$. Additionally, define $V:(\R^+)^2 \to \R$ by $V(x,y)=\frac12 (x^2+y^2)$ and let $(C(t),R(t),P_C(t),P_R(t))$ be a solution of~\eqref{eq:modelo} with initial conditions $C(t_0)=C_0$, $R(t_0)=R_0$, $P_R(0)=P_{R,0}$ and $P_C(0)=P_{C,0}$ and $(C_a(t),R_a(t))$ be a solution of~\eqref{eq:sistema_equiv_para_C_and_R_=_0} with initial conditions $C_a(t_0)=C_0$ and $R_a(t_0)=R_0$. Define $x(t)=C(t)-C_a(t)$ and $y(t)=R(t)-R_a(t)$.

By~\eqref{eq:C_and_R} and~\eqref{eq:sistema_equiv_para_C_and_R_=_0} we have
\begin{equation}\label{eq:estim_x}
\begin{split}
x'
& =C'-C_A'\\
& =\lbd_7 (R-R_a)-(\eps+\beta_1+\lbd_5)(C-C_a)+(\lbd_1+m\lbd_4)(\theta_2(t)-(C-C_a))\\
& \quad + q\lbd_2(R-R_a) +\lbd_2\theta_2(t)R+\lbd_2(R_aC_a-RC)\\
& =\lbd_7 y-(\eps+\beta_1+\lbd_5)x+(\lbd_1+m\lbd_4)(\theta_2(t)-x) + q\lbd_2y +\lbd_2\theta_2(t)R\\
& \quad +\lbd_2(R_a(C_a-C)+(R_a-R)C)\\
& = -(\eps+\beta_1+\lbd_5+\lbd_1+m\lbd_4+\lbd_2R_a)x+ (\lbd_7+q\lbd_2-C\lbd_2)y\\
& \quad + (\lbd_2R+(\lbd_1-m\lbd_4))\theta_2(t)
\end{split}
\end{equation}
and
\begin{equation}\label{eq:estim_y}
\begin{split}
y'
& =R'-R_a'\\
& =\lbd_5 (C-C_a)-(\eps+\beta_2+\lbd_7+\lbd_3+m\lbd_4+(\lbd_2+m_R\lbd_6)p)(R-R_a)\\
& \quad +(\lbd_2+m_R\lbd_6)(R^2-R_a^2)+((\lbd_2+m_R\lbd_6)R+\lbd_3+m\lbd_4)\theta_1(t)\\
& =\lbd_5 x-(\eps+\beta_2+\lbd_7+\lbd_3+m\lbd_4+(\lbd_2+m_R\lbd_6)p)y\\
& \quad +(\lbd_2+m_R\lbd_6)(R+R_a)y+((\lbd_2+m_R\lbd_6)R+\lbd_3+m\lbd_4)\theta_1(t)\\
\end{split}.
\end{equation}
By~\eqref{eq:elim_Pr_Pc} in the proof of Theorem~\ref{teo:general_system}, there is $T_0 \ge 0$ such that $P_R(t),R(t) \le P_R(t)+R(t) < p+\delta$ and
$P_C(t),C(t) \le P_C(t)+C(t) < q+\delta$, for all $t \ge T_0$.
Using this fact, by~\eqref{eq:estim_x} and~\eqref{eq:estim_y} we have
\[
\begin{split}
\frac{d}{dt} V(x,y)
& = xx'+yy'\\
& = -(\eps+\beta_1+\lbd_5+\lbd_1+m\lbd_4+\lbd_2R_a)x^2+
        (\lbd_7+q\lbd_2-C\lbd_2)xy\\
&   \quad + (\lbd_2R+(\lbd_1-m\lbd_4))x\theta_2(t) + \lbd_5 xy\\
&   \quad -(\eps+\beta_2+\lbd_7+\lbd_3+m\lbd_4+(\lbd_2+m_R\lbd_6)p)y^2\\
&   \quad +(\lbd_2+m_R\lbd_6)(R+R_a)y^2+
        ((\lbd_2+m_R\lbd_6)R+\lbd_3+m\lbd_4)y\theta_1(t) \\
& \le -(\eps+\beta_1+\lbd_5+\lbd_1+m\lbd_4)x^2+(\lbd_5+\lbd_7+q\lbd_2-C\lbd_2)xy\\
&   \quad -(\eps+\beta_2+\lbd_7+\lbd_3+m\lbd_4+
    (\lbd_2+m_R\lbd_6)(p+2\delta))y^2+\Theta(t),
\end{split}
\]
for $t \ge T_0$, where
\small{$$\Theta(t)=(\lbd_2R+(\lbd_1-m\lbd_4))(p+\delta)\theta_2(t)+ ((\lbd_2+m_R\lbd_6)R+\lbd_3+m\lbd_4)(q+\delta)\theta_1(t).$$}

Thus, using the fact that $xy \le 1/2(x^2+y^2)$, we get
\[
\begin{split}
\frac{d}{dt} V(x,y)
& \le -(\eps+\beta_1+\lbd_5+\lbd_1+m\lbd_4)x^2+
    \frac12(\lbd_5+\lbd_7+q\lbd_2)(x^2+y^2)\\
&   \quad -(\eps+\beta_2+\lbd_7+\lbd_3+m\lbd_4+
    (\lbd_2+m_R\lbd_6)(p+2\delta))y^2+\Theta(t)\\
& \le -(\eps+\beta_1+\frac12(\lbd_5-\lbd_7)+\lbd_1+m\lbd_4-\frac12q\lbd_2)x^2\\
&   \quad -(\eps+\beta_2+\frac12(\lbd_7-\lbd_5)+\lbd_3+m\lbd_4+
    (\lbd_2+m_R\lbd_6)(p+2\delta)\\
& \quad -\frac12q\lbd_2)y^2+\Theta(t),
\end{split}
\]
for $t \ge T_0$. By~\eqref{eq:cond-equiv}, there is $M > 0$ such that
\[
\frac{d}{dt} V(x,y) \le -M V(x,y) +\Theta(t).
\]

Since $\Theta(t)\to 0$ as $t\to 0$, there is $T_1 \ge T_0 > 0$ sufficiently large such that, for $t \ge T_1$, we have $\Theta(t) \le \delta$. Thus
$$ \frac{d}{dt} V(x(t),y(t)) \le -M V(x,y) + \delta$$
for $t\ge T_1$. Therefore, for $t \ge T_1$,
$$ V(x(t),y(t)) \le \frac{\delta}{M} + V(x(T_1),y(T_1)) \e^{-M(t-T_1)}.$$

We conclude that
 $$\lim_{t \to +\infty} V(x(t),y(t)) \le \frac{\delta}{M}.$$
Since $\delta>0$ is arbitrary, the theorem follows.
\subsection{Proof of Theorem~\ref{teo:static_social_network}}
Adding equations for $R$ and $P_R$ and for $C$ and $P_C$ in~\eqref{eq:modelo} we get the system
\[
\begin{cases}
(R+P_R)'=\alpha\gamma-\eps(R+P_R)\\
(C+P_C)'=(1-\alpha)\gamma-\eps(C+P_C)
\end{cases}
\]
and thus
\[
\begin{cases}
R(t)+P_R(t)=\frac{\alpha\gamma}{\eps}+C_1\e^{-\eps t}\\[2mm]
C(t)+P_C(t)=\frac{(1-\alpha)\gamma}{\eps}+C_2\e^{-\eps t}
\end{cases}.
\]
Therefore, if $(R^*,C^*,P_R^*,P_C^*)$ is an equilibrium solution, then \begin{equation}\label{eq:PR_PC_equilib_static}
    P_R^*=\alpha\gamma/\eps-R^* \quad \quad \text{and} \quad \quad P_C^*=(1-\alpha)\gamma/\eps-R^*.
\end{equation}
Using these expressions and the first two equations in System~\eqref{eq:modelo}, we obtain
\begin{equation}\label{eq_equilibrium_static}
\begin{cases}
-(\eps+\beta_1) C^* + (\lbd_1+m\lbd_4+\lbd_2R^*) \left(\frac{(1-\alpha)\gamma}{\eps}- C^* \right)=0\\
-(\eps+\beta_2)R^* + (\lbd_3+m\lbd_4+(\lbd_2+m_R\lbd_6)R^*) \left(\frac{\alpha\gamma}{\eps}- R^*\right)=0
\end{cases},
\end{equation}
and using the second equation in~\eqref{eq_equilibrium_static} and dividing by $-(\lbd_2+m_R\lbd_6)$ (that is nonzero by assumption), we get
$$(R^*)^2+2\theta R^*-\frac{\alpha\gamma(\lbd_3+m\lbd_4)}{\eps(\lbd_2+m_R\lbd_6)}=0.$$
Thus we have a unique nonnegative root given by~\eqref{eq:static-R*}.
By the first equation in~\eqref{eq_equilibrium_static} we obtain
$$C^* = \frac{(1-\alpha)\gamma}{\eps} \frac{\lbd_1+m\lbd_4+\lbd_2R^* }{\eps+\beta_1+\lbd_1+m\lbd_4+\lbd_2R^*},$$
by~\eqref{eq:PR_PC_equilib_static} we conclude that
$$P_R^*=\dfrac{\alpha\gamma}{2\eps}+\theta-\sqrt{\left(\dfrac{\alpha\gamma}{2\eps}+\theta\right)^2
-\frac{\alpha\gamma(\lbd_3+m\lbd_4)}{\eps(\lbd_2+m_R\lbd_6)}}$$
and it is easy to check that the expression above is nonnegative. Again by~\eqref{eq:PR_PC_equilib_static} we get $P_C^*=\gamma(1-\alpha)/\eps-C^*$ and it is also immediate that this expression is nonnegative.

To study the stability of the equilibrium, we consider the Jacobian matrix $J$ at the equilibrium. Namely
$$
J=\left[
\begin{array}{cccc}
-\eps-\beta_1 & \lbd_2 P_C^* & \lbd_1+m\lbd_4+\lbd_2 R^* & 0 \\
0 & -\eps-A & 0 &  B \\
\beta_1 & -\lbd_2 P_C^* & -\eps-\lbd_1-m\lbd_4-\lbd_2 R^* & 0 \\
0 & A & 0 &-\eps-B\\
\end{array}
\right],
$$
where $A=\beta_2-(\lbd_2+m_R\lbd_6)P_R^*$ and $B=\lbd_3+m\lbd_4+(\lbd_2+m_R\lbd_6)R^*$. It is easy check that the eigenvalues of $J$ are $-\eps$, $-(A+B+\eps)$ and $-(\beta_1+\eps+\lbd_1+\lbd_4m-\lbd_2R^*)$. Thus, all eigenvalues have negative real part and we conclude that the equilibrium solution is locally asymptotically stable.
\subsection{Proof of Theorem~\ref{teo:mri}}
By the second equation in~\eqref{eq:modelo} we can see that, if $(R^*,C^*,P_R^*,P_C^*)$ is an equilibrium solution of our system, then
    $$((\lbd_2+m_R\lbd_6)P_R^*-\eps-\beta_1)R^*=0$$
and thus $R^*=0$ or $P_R^*=(\eps+\beta)/(\lbd_2+m_R\lbd_6)$. Using the remaining equations and condition~\eqref{eq:tau} it is straightforward to obtain the equilibrium points.

In the equilibrium $(R_1^*,C_1^*,P_{R,1}^*,P_{C,1}^*)$, the Jacobian matrix $J$ is given by
$$
J=\left[
\begin{array}{cccc}
-\eps-\beta_1 & \lbd_2 (1-\alpha) \gamma/ \eps &  0 & 0\\
0 & -\eps-\beta_2+(\lbd_2+m_R\lbd_6) \gamma\alpha / \eps & 0 & 0\\
\beta_1 & -\lbd_2 (1-\alpha)\gamma / \eps &  -\eps & 0\\
0 & \beta_2-(\lbd_2+m_R\lbd_6) \gamma\alpha / \eps &  0 & -\eps \\
\end{array}
\right].
$$
We can easily check that the constants $-\eps$, $-b_1-\eps$ and
$\gamma\alpha(\lbd_2+\lbd_6m_R) /\eps - (\beta_2+\eps)$
are the eigenvalues of $J$. Therefore if $\tau>1$ the equilibrium is unstable and if $\tau <1$ it is locally asymptotically stable.

In the equilibrium $(R_2^*,C_2^*,P_{R,2}^*,P_{C,2}^*)$, the Jacobian matrix $J$ has the following form
$$
J=\left[
\begin{array}{cccc}
-\eps-\beta_1 & \lbd_2 P_{C,2}^*  & \lbd_2 R_2^* &  0 \\
0 & 0 & 0 & (\lbd_2+m_R\lbd_6)R_2^* \\
\beta_1 & -\lbd_2 P_{C,2}^* &  -\eps-\lbd_2 R_2^* & 0\\
0 & -\eps & 0 & -\eps-(\lbd_2+m_R\lbd_6)R_2^* \\
\end{array}
\right].
$$
We can easily check that the negative constants $-\eps$, $-(\beta_1+\eps+\lambda_2R_2^*)$ and $-(\lbd_2+m_R\lbd_6)R_2^*$ are the eigenvalues of $J$. Therefore the equilibrium is locally asymptotically stable.
\subsection{Proof of Theorem~\ref{teo:no_referral_influence}}
In our case, equation~\eqref{eq:solution-R} in Theorem~\ref{teo:equilibriums} has a unique solution which is nonnegative and is given by
$$R^*=\dfrac{(\lbd_3+m\lbd_4)up+(\lbd_1+m\lbd_4)\lbd_5q}{uv-\lbd_5\lbd_7}.$$
Thus, again by Theorem~\ref{teo:equilibriums}, the possible equilibrium solutions are given by $(R^*,C^*,P_R^*,P_C^*)$ where $R^*$ is the constant above,
\[
\begin{split}
C^*
&=\dfrac{(\eps+\beta_2+\lbd_7+\lbd_3+m\lbd_4)R^*-(\lbd_3+m\lbd_4)p}{\lbd_5}\\
&=\dfrac{uvp(\lbd_3+m\lbd_4)+\lbd_5vq(\lbd_1+m\lbd_4) - uvp(\lbd_3+m\lbd_4)+\lbd_5\lbd_7p(\lbd_3+m\lbd_4)}{\lbd_5(uv-\lbd_5\lbd_7)}\\
&=\dfrac{(\lbd_1+m\lbd_4)vq + (\lbd_3+m\lbd_4)\lbd_7p}{uv-\lbd_5\lbd_7}.
\end{split}
\]
\[
P_R^*=p-R^*
=\dfrac{up(\eps+\beta_2+\lbd_7)-\lbd_5(\lbd_7p+(\lbd_1+m\lbd_4)q)}
{uv-\lbd_5\lbd_7}
\]
and
\[
P_C^*=p-C^*
=\dfrac{vq(\eps+\beta_1+\lbd_5)-\lbd_7(\lbd_5q+(\lbd_3+m\lbd_4)p)}
{uv-\lbd_5\lbd_7}.
\]
Since $\max\{\kappa_1,\kappa_2\}\le 1$, is immediate that $P_R^* \ge 0$ and $P_C^* \ge 0$. We conclude that there is a unique equilibrium solution $(R^*,C^*,P_R^*,P_C^*)$.

The Jacobian matrix $J$ at the equilibrium is given by
$$
J=\left[
\begin{array}{cccc}
 -\eps-\beta_1-\lbd_5 & \lbd_7 & \lbd_1+m\lbd_4 & 0 \\
\lbd_5 & -\eps-\beta_2-\lbd_7 & 0 &  \lbd_3+m\lbd_4 \\
\beta_1 & 0 & -\eps-\lbd_5-\lbd_1-m\lbd_4 & \lbd_7 \\
0 & \beta_2 & \lbd_5 & -\eps-\lbd_7-\lbd_3-m\lbd_4 \\
\end{array}
\right].
$$
We can check that the eigenvalues are all negative and given by $-\eps$, $-\eps-\lbd_5-\lbd_7$ and
    $$\frac12 \left( -\sigma -2\lbd_4 m \pm \sqrt{\sigma^2-4[(\eps+\beta_1+\lbd_1+\lbd_5)(\eps+\beta_2+\lbd_3+\lbd_7)-\lbd_7\lbd_5]} \right).$$
Therefore the equilibrium is asymptotically stable. Since the system in this case is linear, the equilibrium is globally asymptotically stable.
\section{Conclusions} \label{section:Conc}

We presented and studied in this paper a compartmental model with four compartments to describe the evolution of the number of customers and potential customers of some corporation based on the marketing policy of the corporation, determined by the effort used in undifferentiated marketing campaigns and in referral directed marketing. Apparently the results obtained are reasonable in the sense that the qualitative behavior obtained is not going against common sense. The results show that the model works in theory in the several scenarios considered, with and without marketing incentives, thus relying on normal marketing policies or on incentives to referrals. The model shows that, in theory, it is possible to predict the influence referrals can have on their peers based on the incentives given to them by the company. The different scenarios also allows to see, in a specific period of time, what happens to the number of current and potential customers based on the zero incentives policy, and therefore, only based on the natural attractiveness power of the referrals (Theorem~\ref{teo:mri} with $\lambda_6=0$); what happens if the company invests in marketing but not on incentives to referrals (Theorem~\ref{teo:no_referral_influence}); and also, what happens if the company only invests on incentives to referrals (Theorem~\ref{teo:mri} with $\lambda_6 \ne 0$). Based on this results, this model allows companies to adjust their marketing policies in order to maximize a specific parameter of the model. For instance, the model allows a company to estimate the amount of investment necessary to transform an $x$ number of potential customers in customers.

This study must be followed by work where the model is used in real world problems. In fact, comparing the results given by the model and real data would of fundamental importance to confirm the usefulness of it. Naturally this is a major task whose feasibility will certainly depend on the accuracy in the estimation of parameters.

 We believe that this work opens several possibilities for future studies. For instance, it would be interesting to consider versions of this model with time-varying parameters to model seasonal phenomena that may occur in some economic activities. It can also be of interest to consider age structured populations, to distinguish different consumption habits, or to subdivide the universe of referrals, to reflect different aspects that make those customers important to the corporation.

\bibliographystyle{elsart-num-sort}

\end{document}